\documentclass[11pt]{article}
\usepackage{amsmath,amssymb,amsfonts}
\def\hybrid{\topmargin 0pt \oddsidemargin 0pt
\headheight 0pt \headsep 0pt

\textwidth 15.4cm \textheight 21.2cm \voffset=0.5cm \hoffset=0.2cm

\marginparwidth 0.0in
\parskip 5pt plus 1pt \jot = 1.5ex}
\catcode`\@=11
\def\marginnote#1{}

\newcommand{\CC}{\ensuremath{\mathbb C}}
\newcommand{\CP}{\ensuremath{\mathbb CP}}

\newcommand{\RR}{\ensuremath{\mathbb R}}

\newcommand{\ZZ}{\ensuremath{\mathbb Z}}

\def\Vsquare#1{\vbox{\Square{$#1$}}\kern-\Thickness}

\def\numberbysection{\@addtoreset{equation}{section}
 \def\theequation{\thesection.\arabic{equation}}}
\numberbysection

\renewcommand{\theequation}{\thesection.\arabic{equation}}
\newcommand{\l@qq}[2]{\addvspace{2em}
 \hbox to\textwidth{\hspace{1em}\bf #1 \dotfill #2}}

\newcounter{app}

\def\app{\setcounter{equation}{0}
\def\theequation{\Alph{app}.\arabic{equation}}\par
 \addvspace{4ex}
 \@afterindentfalse
 \secdef\@app\@dapp}

\newcommand\@app{\@startsection {app}{1}{0ex}%
 {-3.5ex \@plus -1ex \@minus -.2ex}%
 {2.3ex \@plus.2ex}%
 {\normalfont\Large\bf}}
\def\@dapp#1{%
{\parindent \z@ \raggedright \bf #1}\par\nobreak}
\def\l@app#1#2{\ifnum \c@tocdepth >\z@
 \addpenalty\@secpenalty
 \addvspace{1.0em \@plus\p@}%
 \setlength\@tempdima{8em}%
 \begingroup
 \parindent \z@ \rightskip \@pnumwidth
 \parfillskip -\@pnumwidth
 \leavevmode \bfseries
 \advance\leftskip\@tempdima
 \hskip -\leftskip
 #1\nobreak\hfil \nobreak\hb@xt@\@pnumwidth{\hss #2}\par
 \endgroup\fi}
\newcounter{sapp}[app]

\def\sapp{\def\theequation{\Alph{app}.\arabic{equation}}
\par
\@afterindentfalse
 \secdef\@sapp\@dsapp}
\newcommand{\@sapp}{\@startsection{sapp}{2}{\z@}%
 {-3.25ex\@plus -1ex \@minus -.2ex}%
 {1.5ex \@plus .2ex}%
 {\normalfont\large\bfseries}}

\def\@dsapp#1{%
{\parindent \z@ \raggedright \bf #1
}\par\nobreak}
\newcommand{\l@sapp}{\@dottedtocline{2}{1.5em}{2.3em}}

\def\titlepage{\@restonecolfalse\if@twocolumn\@restonecoltrue\onecolumn
 \else \newpage \fi \thispagestyle{empty}\c@page\z@
 \def\thefootnote{\fnsymbol{footnote}} }

\def\endtitlepage{\if@restonecol\twocolumn \else \fi
 \def\thefootnote{\arabic{footnote}}
 \setcounter{footnote}{0}} 
\relax

\hybrid

\newtoks\@stequation

\def\subequations{\refstepcounter{equation}%
 \edef\@savedequation{\the\c@equation}%
 \@stequation=\expandafter{\theequation}
 \edef\@savedtheequation{\the\@stequation}
 \edef\oldtheequation{\theequation}%
 \setcounter{equation}{0}%
 \def\theequation{\oldtheequation\alph{equation}}}

\def\endsubequations{%
 \setcounter{equation}{\@savedequation}%
 \@stequation=\expandafter{\@savedtheequation}%
 \edef\theequation{\the\@stequation}%
 \global\@ignoretrue}

\makeatletter
\newdimen\normalarrayskip 
\newdimen\minarrayskip 
\normalarrayskip\baselineskip
\minarrayskip\jot
\newif\ifold \oldtrue 
\def\arraymode{\ifold\relax\else\displaystyle\fi} 
\def\eqnumphantom{\phantom{(\theequation)}} 
\def\@arrayskip{\ifold\baselineskip\z@\lineskip\z@
 \else
 \baselineskip\minarrayskip\lineskip1\baselineskip\fi}

\def\@arrayclassz{\ifcase \@lastchclass \@acolampacol \or
\@ampacol \or \or \or \@addamp \or
 \@acolampacol \or \@firstampfalse \@acol \fi
\edef\@preamble{\@preamble
 \ifcase \@chnum
 \hfil$\relax\arraymode\@sharp$\hfil
 \or $\relax\arraymode\@sharp$\hfil
 \or \hfil$\relax\arraymode\@sharp$\fi}}

\def\@array[#1]#2{\setbox\@arstrutbox=\hbox{\vrule
 height\arraystretch \ht\strutbox
 depth\arraystretch \dp\strutbox
 width\z@}\@mkpream{#2}\edef\@preamble{\halign \noexpand\@halignto
\bgroup \tabskip\z@ \@arstrut \@preamble \tabskip\z@ \cr}%
\let\@startpbox\@@startpbox \let\@endpbox\@@endpbox
 \if #1t\vtop \else \if#1b\vbox \else \vcenter \fi\fi
 \bgroup \let\par\relax
 \let\@sharp##\let\protect\relax
 \@arrayskip\@preamble}
\def\eqnarray{\stepcounter{equation}%
 \let\@currentlabel=\theequation
 \global\@eqnswtrue
 \global\@eqcnt\z@
 \tabskip\@centering 
 \let\\=\@eqncr
 $$%
 \halign to \displaywidth \bgroup
 \eqnumphantom \@eqnsel
 \hskip\@centering 
 $\displaystyle \tabskip\z@ {##}$%
 &\global\@eqcnt\@ne \hskip 2\arraycolsep
 $ \displaystyle \arraymode{##}$\hfil
 &\global\@eqcnt\tw@ \hskip 2\arraycolsep
 $\displaystyle\tabskip\z@{##}$\hfil
 \tabskip\@centering
 &{##}\tabskip\z@\cr}
\makeatother
\newtheorem{theorem}{Theorem}[section]
\newtheorem{definition}{Definition}[section]
\newtheorem{proposition}{Proposition}[section]

\newtheorem{remark}{Remark}[section]

\def\al{\alpha}
\def\la{\lambda}

\def\<{\langle}
\def\>{\rangle}

\def\g{{\mathfrak g}}
\def\a{{\mathfrak a}}
 \def\hg{\hat{\mathfrak g}}

\def\h{{\mathfrak h}}

\def\n{{\mathfrak n}}
\def\b{{\mathfrak b}}

\def\sln{\mathfrak{sl}_n}
\def\slt{\mathfrak{sl}_2}
\def\slf{\mathfrak{sl}_3}
\def\hgv{\mathfrak{ g}^V}

\def\cd{{\cal D}}
\def\ca{{\cal A}}

\def\1{1\!\!1}

\def\tp{\otimes}

\def\d{\partial}

\begin{document}
\thispagestyle{empty}

\begin{titlepage}
 \bigskip\bigskip
 \begin{center}
 \bigskip\bigskip
 {\Large\bf Higher-Dimensional generalizations of Affine Kac-Moody and
 Virasoro Lie Algebras}\\
 \bigskip\bigskip
 {\large M.Golenishcheva-Kutuzova \footnote{E-mail: kutuzova@math.ufl.edu}}
 \\ \bigskip
 \textit{Department of Mathematics, University of Florida \\ PO Box 118105,
 Gainesville FL 32611-8105}\\
 \bigskip
 \end{center}
 \bigskip \bigskip \bigskip
 \begin{abstract}
 We discuss the higher dimensional generalizations of the Virasoro and the Affine
 Kac-Moody Lie algebras. We present an explicit construction for a central
 extensions of the Lie Algebra $Map\,(X, \g)$ where $\g$ is a
 finite-dimensional Lie algebra and $X$ is a complex manifold that can be
 described as a "right" higher-dimensional generalization of $\CC^*$ from the
 point of view of a corresponding group action. The constructed algebras have
 most of the good properties of finite dimensional semi-simple Lie algebras and
 are a new class of generalized Kac-Moody algebras. These algebras have
 description in terms of higher dimensional local fields.

 \end{abstract}
\end{titlepage}
\tableofcontents
\section{Introduction}

The basic object of $d$-dimensional Quantum Field Theory (QFT) are
$d$-dimensional fields (operator-valued distributions on $M^d$)
and the representation of the ${\rm Poincar\acute e}$ group in the
Hilbert Space $H$ that behave accordingly to the Wightman axioms
of QFT. In 1984 Belavin, Polyakov and Zamolodchikov (\cite{BPZ})
initiated the study of two-dimensional Conformal Field Theory
(CFT). In CFT we have a nice split of variables that leads to the
notion of the "chiral part" of a conformal field theory, where
fields are just operator-valued distributions on $\CC^*.$ This
rises the huge interest to the theory of infinite-dimensional Lie
algebras and their representations. It turns out that the
symmetries of the most models of CFT are classified along the
representations of Affine Kac-Moody Lie algebras and the Virasoro
algebra.

A mathematical definition of the "chiral part" of a conformal
field theory, called a Vertex algebra, was proposed by Borcherds
(\cite{B}). The axioms of Vertex algebras (\cite{K2}) are
mathematical description of the operator product expansion in CFT.
The Vertex algebra approach simplifies the problem of
classification of infinite dimensional symmetries in CFT. However,
until now a classification of Vertex algebras seems to be far
away. There is the solution to the classification problem only
when the chiral algebra is generated by a finite number of quantum
fields, closed under the operator product expansions (in a sense
that only derivatives of the generating field may occur). Roughly
speaking, in the finitely generated case there is noting but the
Affine Kac-Moody algebras and the Virasoro Vertex algebra
(Heisenberg algebra can be treated as a trivial Affine algebra).
This fact was proven by Kac using the notion of a conformal
algebra (\cite{K2}), which is related to a chiral algebra in the
same way as a Lie algebra is related to its universal enveloping
algebra. On the space of fields (mutually local formal
distributions in complex variable $z$) there is the action of the
operator $\d,$ given by $\d a(z)=\d_z a(z)$ and the fields that
appear in the right part of the commutation relations of two
fields $a(z)$ and $b(w)$ can be viewed as $n$-products of these
two fields.

The definition of conforormal algebra is:
 \begin{definition}
 [V.Kac]
 A (Lie) conformal algebra is a $\CC[\d]$-module $R$, endowed with a family of
 $\CC$-bilinear products $a_{(n)}b, \,\, n\in \ZZ_+,$ satisfying axioms
 (C1)--(C4):

 (C1) $a_{(n)}b =0 \,\,\, for \,\,\, n>>0,$

 (C2) $(\d a)_{(n)}b = -n a_{(n-1)}b, \,\,\, a_{(n)}\d b= \d (a_{(n)}b) - (\d
 a)_{(n)}b,$

 (C3) $a_{(n)}b=- \sum_{j=0}^{\infty} \, (-1)^{n+j} \d^{(j)}(b_{(n+j)}a),$

 (C4) $ a_{(m)}(b_{(n)}c)- b_{(n)}(a_{(m)}c) = \sum_{j=0}^{m} \,
 {m \choose j} (a_{(j)}b)_{(m+n-j)} c,$\\
 where $\d^{(j)}=\frac{ \d^{j}}{n!}.$

 A conformal algebra $R$ is called $finite$ if $R$ is a finitely generated
 $\CC[\d]$-module. The $rank$ of conformal algebra $R$ is its rank as a
 $\CC[\d]$-module.
 \end{definition}
 \begin{definition} Let $\g$ be an arbitrary Lie algebra. A formal
 distribution Lie algebra $(\g, \, F)$ is the space $F$ of all
 mutually local $\g$-valued distributions in complex  variable $z.$
\end{definition}

It is very important that we have a functor from the category of
formal distribution Lie algebras to the category of conformal
algebras as well as a functor in the opposite direction that
canonically associates to a conformal algebra $R$ a formal
distribution Lie algebra (\cite{K2}). It means that when we use a
very formal language of conformal algebras we do not loose the
information about the physical origin. It is specially important
to have it in mind when we would like to construct higher
dimensional generalizations.

The approach to higher dimensional chiral algebras suggested by
Beilinson and Drinfeld is based more on algebraic geometry than
the representation theory. In (\cite{BD}) they introduced the
notion of "Chiral algebra" as a quantization of what they call the
"coisson algebra" (a Poisson algebra on $X$ in the compound
setting). A really challenging problem is to find out what is a
chiral algebra on higher dimensional $X$ (the coisson algebras
live in any dimension). More algebraic approach to the higher
dimensional conformal algebras was suggested by Bakalov, D'Andrea,
and Kac in (\cite{BDK}, \cite{BDK1}) introducing the notion of
$Lie \: \: pseudoalgebra.$ The basic idea is to replace the
$\CC[\d]$ in the definition of conformal algebra $R$ by a Hopf
algebra $H = U(\g),$ where $\g$ is a finite-dimensional Lie
algebra and $U(\g)$ is its universal enveloping algebra. A Lie
pseudoalgebra is defined as an $H$-module $L$ endowed with an
$H$-bilinear map
\begin{equation}
 L\tp L \to (H\tp H)_{\tp_H } \,L
\end{equation}
subject to a certain skewsymmetry and Jacobi identity axioms.
Unfortunately, for Lie pseudoalgebras the functor to the higher
dimensional distribution Lie algebras was not explicitly
constructed.

We suggest an approach to the higher dimensional CFT that
preserves the connection between the higher dimensional
distribution Lie algebras and a higher dimensional conformal
algebras. The right idea is to replace the $\CC[\d]$ in the
definition of conformal algebra $R$ by a universal enveloping
algebra of some noncommutative Lie algebra $\n.$ We consider the
case when  $\n$ is a nilpotent subalgebra of a simple Lie algebra
$\g.$ In this situation we can construct a Lie algebra of formal
local distributions on some complex manifold $M^{\g}$ that has the
most of the properties of two dimensional CFT. It means that we
have the consistent definition of the OPE and the normal product
of two fields in the dimension higher than two. The corresponding
conformal algebra is a $U(\n)$-module with the
$e_{\gamma}$-products subject to a certain skewsymmetry and Jacobi
identity axioms. Here $\{e_{\gamma}\}$ is a basis in some
representation of $\g.$ When $\g=\slt$ our construction leads to
the standard two dimensional situation, where $M^{\g}= \CC^*$ and
the $n$-product of two fields corresponds to the $(e_n)$-product,
where $\{e_n\}$ is a basis in some representation of $\slt$.

As it was mentioned above, in the case of a finite conformal
algebra, there are only two principal solutions to the system of
axioms (C1)--(C4) (in this paper we do not discuss a conformal
superalgebras, where we have more possibilities):

1. {\bf Current conformal algebras Cur $\g =\CC [\d] \tp \g$
associated to the Lie algebra $\g.$} The only non-trivial
$n$-product is the $0$-product: $a_{(0)}b =[a,\, b], \, \, a, b
\in \g.$ We identify $\g$ with the subspace of Cur $\g$ spanned by
elements $1 \tp g, \, \, g \in \g.$ The corresponding formal
distribution Lie algebra is $(\tilde{\g}, R)$ where $\tilde{\g}=\g
\tp \CC[t, \,t^{-1}].$   Formal distributions $g(z)=g \tp \delta
(t-z)=\sum_{n\in\ZZ} \, g\tp t^n \tp z^{-n-1},$ defined for every
$g\in \g$ and satisfy the commutation relations:
\begin{equation}
 [g_1(z), \, g_2(w)]= [g_1,\, g_2](w) \delta (z-w).
\end{equation}

2. {\bf Virasoro conformal algebra Conf (Vect $\CC^*).$ } The
centerless Virasoro algebra of algebraic vector fields on $\CC^*$
is spanned by the vector fields $L_n = -t^{n+1}\, \d_{t}.$ The
Vect$\CC^*$-valued formal distribution $L(z) = \delta(t-z)\,
\d_{t}$ satisfies
\begin{equation}
 [L(z), \, L(w)] = \d_w L(w) \,\delta(z-w) + 2 L(w) \, \d_w \delta (z-w).
\end{equation}
The conformal algebra Conf ${\rm(Vect}\, \CC^*)= \CC[\d] L$
associated to (Vect $\CC^*, \{L(z)\})$ is defined by the
$n$-products:
\begin{equation}
 \label{vir}
 L_{(0)}L= \d L, \,\,\,\,\;\; L_{(1)}L= 2L, \,\,\,\;\;\; L_{(n)}L= 0 \,\,\,\;{\rm
 if}\,\,\,\;\; n \ge 2.
\end{equation}

The conformal algebra associated to the semidirect sum (Vect$\, \CC^*) +
\tilde{\g}$ is the semidirect sum $\,$ Conf(Vect $\CC^*)$ + Cur $\g,$ defined by
$(a \in \g):$
\begin{equation}
 L_{(0)}a = \d a, \,\,\,\, L_{(1)}a= a, \,\,\,\, L_{(n)}a=0 \,\,\, {\rm
 for}\,\,\,\, n>1.
\end{equation}

What makes these two cases so special? The answer is that they
both are related to the action of the Lie algebra $\slt$ on the
space $V^{\slt}$ of regular functions on a subset $M\simeq \CC^*$
of the flag manifold $SL(2,\CC)/B_- \simeq \CP^1$ and the
corresponding fields are expressed in terms of the formal delta
function $\delta(t-w)= \sum_{n\in \ZZ}\;{t^n}\,{w^{-n-1}},$
associated with this space.

Our approach is based on the construction of a complex manifold
$M^{\g}$ that can be described as a "right" higher-dimensional
generalization of $\CC^*$ from the point of view of a
corresponding group action. In  sect.2 we define the space
$V^{\g}$ for any simple Lie algebra $\g$ and the formal delta
function associated with this space. In sect.3 we define the Lie
algebra of formal distributions on $V^\g$ and the higher
dimensional analogues of the $n$-products and $\lambda$-product
for the conformal algebra associated with the Lie algebra of
formal distributions on $V^\g$.

In sect.4  we discuss the basic ideas about the geometrical
realization
 of the space $ V^{\g}$ in general case for any
simple Lie algebra $\g ,$  and in sect.5 we present an explicit
realization of  $ V^{\g}$ for $\g=\slf .$

In sect.6 we define  the  3-dimensional analogues of the Affine
Kac-Moody algebras associated with the action of the Lie algebra
$\slf .$ We call these algebras the Generalized Affine Kac-Moody
algebras $\hgv$ associated with the space $V^{\slf}.$ We think
that it is especially important to consider with more details the
$\slf$-case, because it is well known fact, that it is difficult
to make a step from $\slt$ to $\slf$ and it is almost strait
forward from $\slf$ to go  to the general case of of any simple
Lie algebra (at least to $\sln$). The constructed algebras have
many of the good properties of "the generalized Affine Kac-Moody
algebras". The basic idea of Borcherds is to think of generalized
Affine Kac-Moody algebras as infinite dimensional Lie algebras
which have most of the good properties of finite dimensional
reductive Lie algebras (\cite{B}). For instance, we construct the
normalized invariant form $(\,,\,)$ and the  Cartan involution of
$\hgv .$

In sect.7 we discuss the higher dimensional version of the
Virasoro conformal algebra associated with the space $V^{\slf}.$
The Virasoro algebra appears in  many different contexts related
to the Lie algebra $\slt$ and contains $\slt$ as a subalgebra. We
will discuss one more context related to the conformal Virasoro
(\ref{vir}) algebras and similarly we define the generalized
Virasoro conformal algebra associated with the space $V^{\slf}.$
This algebra contains $\slf$ as a subalgebra and it is a rank one
module over $U(\n_+)$, where $\n_+$ is the upper nilpotent
subalgebra of $\slf .$ Also we have the semidirect sum of the
generalized Virasoro conformal algebra and conformal algebra $\hgv
$. The generalized Virasoro conformal algebra (the root lattice
Virasoro conformal algebra) associated to any simple Lie algebra
$\g$ will be constructed in (\cite{G-K}). We would like to notice
that our definition of the generalized Virasoro conformal algebras
is different from the Virasoro psuedoalgebras defined in
(\cite{BDK}, \cite{BDK1}).

In the future publications we shall present:\\
1. The explicit geometrical realization of the space $V^\g$ for
any simple Lie algebra $\g$;\\
 2. The representations of the
Generalized Affine Kac-Moody algebras $\hgv$ associated with the
space $V^{\g}$ and the generalized Virasoro conformal algebra
associated with this space;\\
3. The axioms for the higher dimensional Vertex Algebras
associated with the space $V^{\g} ;$\\
4. The higher dimensional generalizations of the $N=1$ and $N=2$
superconformal algebras.

{\bf Acknowledgments.} The author wants to thank A. Gerasimov for
very useful discussions and M. Olshanesky for reading the
manuscript and suggesting improvements.

\section{Generalized delta-function}

In this section we introduce the notion of the generalized
delta-function associated with a vector space $V^{\g},$ where $\g$
is a simple Lie algebra. First we recall the basic definitions of
formal distributions and in particular the formal delta-function.

A formal distribution in the indeterminates $z, w, \ldots \in (\CC^*)^N$ with
values in a vector space $W$ is a formal expression of the form
\begin{equation*}
 \sum_{m,n\ldots , \in \ZZ} a_{m,n,\ldots } \, z^m w^n \ldots ,
\end{equation*}
where $a_{m,n,\ldots }$ are elements of a vector space $W$. They form a vector
space denoted by $W[[z, z^{-1}, w, w^{-1}]]$.

Given a formal distribution $ a(z) =\sum_{m\in \ZZ} \, \, a_m
z^m$, we define the trace (integral) by the usual formula
\begin{equation}\label{tree}
 Res_{z}\, a(z) = a_{-1}.
\end{equation}
Define the $\CC$ - valued bilinear form  on on the space of $\CC$
- valued formal distribution by
\begin{equation}
 \label{tr}
 <a(z),\, b(z)> = Res_{z}\, a(z)b(z)
\end{equation}
Since $Res_{z}\, \d a(z) =0$, we have the usual integration by part:
\begin{equation}
 \label{tr1}
 <\d a(z),\, b(z)> = -<a(z),\, \d b(z)>,
\end{equation}
where $\d a(z)= \frac{\d}{\d z}\,a(z)$.

We would like to remark, that the bilinear form (\ref{tr}) is invariant under
the action of the Lie algebra $\slt$ given by:

\begin{equation}
 \label{rep}
 X= \d ,\;\;\;\;\;\;
 H= -2 z \d -1 ,\;\;\;\;\;
 Y= -z^2\d - z,
\end{equation}
where $X,H,Y$ is the standard basis in $\slt$. This bilinear form
defines the pairing between the Verma module of $\slt$ with the
highest weight $-\frac{\al}{2}$ in the space $V^+ \simeq \CC [z]$
and the Verma module with the lowest weight $\frac{\al}{2}$ in the
space $V^- \simeq \frac{1}{z}\, {\CC [\frac{1}{z}]}.$

The action (\ref{rep}) results from the natural action of the
group $G=SL(2, \CC)$ on the Flag manifold $ G/B_{-}$, where
$B_{-}$ is the lower Borel subgroup of $G .$  The space $V^+$ is
isomorphic to the space of regular functions on the big cell $U=
 N_+ \cdot [1] \in G/B_{-}$ and the space $V^-$ is isomorphic to the
space of regular functions on the dual cell $U^*= N_- \cdot
s_{\al}[1] \in G/B_{-}$ factorized by constants. Here $s_{\al}$
denote the action of the generator of the Weyl group of $\slt$ on
the flag manifold. The space $V^+$ is the maximum isotropic
subspace with respect to the bilinear form (\ref{tr}). The weight
basis in $V^+$ is ${e_{n}=z^n},\, \, \, n \in \ZZ_{+}$ and the
dual basis in $V^{-}$ is ${e_{n}}^*=z^{-n-1}, \, \, n \in
\ZZ_{+}.$ The complex torus $\CC^*$ is the intersection of $U$ and
$U^*$ and the space of all regular functions on $\CC^*$ is $\CC[z,
z^{-}]= V^+ \oplus V^- .$

The subspaces $V^+$ and $V^-$ are invariant under the multiplication and the
action of $\slt$, given by (\ref{rep}). We will denote the space $\CC[z,
z^{-1}]$ by $V_{z}^{\g}$ for $\g = \slt$.

Recall that the formal delta function is the following formal distribution in
$z$ and $w$ with values in $\CC$
\begin{equation}
 \delta (z-w) =\sum_{n\in \ZZ}z^n \cdot w^{-n-1}.
\end{equation}
We can think about the formal delta function as an element of the space
$V_{z}^{\g} \tp V_{w}^{\g}$ of the form
\begin{equation}
 \label{del}
 \delta (z-w) =\sum_{n\in \ZZ_+} e_{n}\tp e^*_{n} + \sum_{n\in \ZZ_+} e^*_{n}\tp
 e^{n}=\\
 \delta (z-w)_- + \delta (z-w)_+,
\end{equation}
where
\begin{equation}
 \label{d+}
 \delta (z-w)_- \in V^+_{z} \tp V_{w}^{-} \, \, \, \rm{and} \,\,\, \delta
 (z-w)_+ \in V^-_{z} \tp V_{w}^{+}.
\end{equation}
The well known properties of the formal delta-function (\cite{K1})
result from  (\ref{del}). Let us mention some of them:

(a) For any formal distribution $f(z) \in U[[z, z^{-1}]]$ one has: $ Res_{z}\,
f(z) \delta (z-w)=f(w),$

(b) $\d^{j}_{z}\, \delta (z-w) =(-\d_{w})^j \, \delta (z-w).$

As we have mentioned before, the formal delta function is
connected with the action of the Lie algebra $\slt$ on functions
on the Flag manifold and it is an element of the space $V_{z}^{\g}
\tp V_{w}^{\g}$ (\ref{rep}).

 The construction of the generalized
formal delta functions is based on the same idea. First we will
give a formal definition of the space $V^{\g}$ for a simple Lie
algebra $\g$ and the generalized formal delta functions connected
with this space. Then in the next sections we will give the
explicit construction for $\g= \slf$ and discuss the general
construction for any simple Lie algebra $\g.$

Let $\g$ be a simple Lie algebra of rank $l$. As a vector space,
it has the triangular decomposition
\begin{equation}
 \g =\n_{+} \oplus \h \oplus \n_-,
\end{equation}
where $\h$ is a Cartan subalgebra and $\n_{\pm}$ are the upper and
lower nilpotent subalgebras. Let
\begin{equation}
\b_{\pm}= \h \oplus \n_{\pm}
\end{equation}
be the upper and lower Borel subalgebras.
\begin{definition}
A space $V^{\g}$ is a vector space endowed with $\CC$-valued
non-degenerated symmetric bilinear form \begin{equation} < \cdot , \, \cdot >
\, : V^{\g} \tp V^{\g} \longrightarrow \CC \end{equation} such that the
following axioms hold

(V1) $V^{\g}$ it a commutative associative algebra with unitary element $\bf{1}$ with respect to
the multiplication. ( We can think that $V^{\g}$ is a space of complex-valued functions on
a complex manifold $M$ or orbifold).

(V2) $V^{\g}$ is a $\g$-module, such that the action of $ \n_+$ is
a derivation of $V^{\g}$; it means that $x(f \cdot g) = x(f)\cdot
g + f \cdot x(g)$ for any $x \in \mathfrak{n}_+$ and any $f, \, g
\in V^{\g}$ . In particulary, it means that the unitary element of
$V^{\g}$ is annihilated by all elements from $ \mathfrak{n}_+$.

(V3) The bilinear form $< \cdot , \, \cdot >$ is invariant under
the multiplication in $V^{\g}$ and the action of the Lie algebra $
\n_+$:
\begin{equation}
<f\cdot g ,\, \, h> = <g , \,\, f\cdot h>\\
<x(f), \, \,g> +< f, \, \,x(g)> =0,
\end{equation}
for any $x \in \mathfrak{n}_+$ and any $f, \, g , h \in V^{\g}$ .

(V4) There is a maximal isotropic with respect to the bilinear form subspace $V_+$ in $V^{\g}$,
such that $V_+$ is invariant under the multiplication and the action of $ \mathfrak{n}_+$.
\end{definition}
Then
\begin{equation}
 V^{\g} = V_ + \oplus V_- ,
\end{equation} where $V_- \simeq V^*_+$ is the dual to
the subspace $V_+$ in $V^{\g}$ with respect to the given bilinear
form.

In the next section we will construct the explicit realization of
the space $V^{\g}$ as the ring of complex-valued regular functions
on a complex manifold $M.$ Because the vector space $V^{\g}$ is
self-dual, we can identify $V^{\g}$ with $ (V^{\g})^*$. Let
$\{e_{\gamma}\},\, \, \, \gamma \in \Gamma$ be a basis in $V_+$
and $\{e_{\gamma}^* \}$ is the dual basis in $V^{-} ,$ where
 $ \Gamma$ is a discrete set that numerate the basis.
\begin{remark}
In the previous definition instead of bilinear form we can
postulate that the space $V^{\g}$ has
 a trace $\, Res$ invariant with respect to the action of the nilpotent subalgebra $\mathfrak{n}_+$.
\end{remark}

 If such
trace exists, the bilinear form on $V^{\g}$ can be defined as
\begin{equation}\label{trace}
 <f \, ,\, \, g> = Res \, (f \cdot g).
 \end{equation}
 In the
other direction, if we have a bilinear form with the given
properties on $V^{\g},$ and ${e_{\gamma_0}}=\bf{1}$ for some
$\gamma_0 \in \Gamma$, then the invariant trace can be defined as
\begin{equation}\label{res}
 Res \, (f )\, = <f \, ,\, \, \bf{1}>= \, \rm{coefficient\ of } \, e_{\gamma_0}^*
\end{equation}
in the decomposition with respect to the basis
$\{e_{\gamma}\}\cup \{e_{\gamma}^*\}, \, \, \gamma \in \Gamma.$

\begin{definition} The generalized formal delta functions associated with the space
 $V^{\g}$ is defined as an element
of the space $V^{\g} \tp (V^{\g})^* = V_+ \tp V_- \oplus V_- \tp V_+$ of the form
\end{definition}
\begin{equation}
 \label{delt}
 \delta_{\scriptstyle{V^{\g}}} =\sum_{\gamma \in \Gamma} e_{\gamma} \tp
 e_{\gamma}^* \, \oplus \, e_{\gamma}^* \tp e_{\gamma}.
\end{equation}

Now suppose that the space $V^{\g}$ is realized as the space
$\rm{Fun}(M)$ of formal distributions on $M$ and ${\bf{z}} =(z_1,
z_2, \ldots , z_n)$ are coordinates on $M,$ then the space $V^{\g}
\tp (V^{\g})^*$ can be identified with the space of distributions
in two sets of coordinates ${\bf{z}} =(z_1, z_2, \ldots , z_n)$
and ${\bf{w}}=(w_1, w_2, \ldots , w_n).$ In this case we will use
the notation
\begin{equation} \delta_{V^{\g}}\,\, (\bf{z}-\bf{w}),
\end{equation} or for short $\delta_{V}\, (\bf{z}-\bf{w}). $

The action of the nilpotent subalgebra $\n_+$ on Fun$(M)$ is given
by vector fields. Thus, we have a Lie algebra homomorphism $\n_+
\to \rm{Vect} M$.  So defined generalized formal delta function
have the most of the properties of the standard delta function
with respect to the trace (\ref{tr}) and differentiations from
$\n_+$:

(a) for any formal distribution $f({\bf z}) \in V^{\g}$ one has:
\begin{equation} Res_{\bf z} \,f({\bf z})\, \delta _V ({\bf z}-{\bf w})=f({\bf
w}), \end{equation}

(b) for any element $a \in \n_+$ one has: \begin{equation} (\xi_a^z)^j \,
\delta ({\bf z}-{\bf w}) =(- \xi_a^{w})^j \, \delta ({\bf z}-{\bf
w}), \end{equation} where $\xi_a^z$ is the image of $a$ in $\rm{Vect}\, M$.

\section{Local distributions on $ V^{\g} $}

 By the definition the space $ V^{\g} = V_ + \oplus V_- $
 is a $U(\n_+)$-module. We will think that $ V^{\g}$ is realized
 as a space of regular functions on some complex manifold $M^{\g}.$
Let $\{e_{\gamma}\},\, \, \, \gamma \in \Gamma$ be a basis in
$V_+$ and let $\{e_{\gamma}^* \}$ be the dual basis in $V_{-} $
constructed in the previous section.

Fix
 a triangular decomposition of $\g$:
\begin{equation}
 \g =\n_{+} \oplus \h \oplus \n_-
\end{equation} where $\h$ is the Cartan subalgebra and $\mathfrak{n}_{\pm}$ are
the upper and the lower nilpotent subalgebras. Let
\begin{equation} \b_{\pm}= \h \oplus \n_{\pm} \end{equation} be
the upper and lower Borel subalgebras. Let $h_i = \al_i, \, \,
i=1,\ldots ,l,$ be the $i$th coroot of $\g$ and let $l$ be the
rank of $\g.$ The set $\{h_i \}_{i=1, \ldots ,l}$ is a basis of
$\h$. We choose a root basis of $\n_{\pm}$ , $\, \{e^{\al}\}_{\al
\in \Delta_{\pm}}$, where $\, \Delta_{+}$ ($\Delta_{-}$) is the
set of positive (negative) roots of $\g$, so that $[ h, \,
e^{\al}] = \al (h) e^\al$ for all $h \in \h.$

Let $ \xi_{\al_1}$ be the image of $e^{\al_i}, \,\,\, \al \in
\Delta_{+}$ in Der$\,(V^\g).$ Fix some ordering of $\Delta_{+} .$
Then we can write $ \d_{i}$ instead of $ \xi_{\al_i}$
 and $ (V^{\g})$ is a $\CC[\d_{1},\ldots ,\d_{p}]$-module, where
$p= |\Delta_+|.$ The $\CC[\d_{1},\ldots ,\d_{p}]$ has a ${\rm
Poincar\acute e}$-Birkhoff-Witt basis of the form
$\{\d_1^{k_1}\,\d_2^{k_2}\ldots \d_p^{k_p}\}.$

 Let $W$ be a vector space. Consider the space ${\rm End}\, W \tp V^{\g}$
of ${\rm End}\, W$-valued formal distributions associated with the
space $V^{\g}.$ Any element $a^{V} \in{\rm End}\, W \tp V^{\g}$ is
an expression of the form
 \begin{equation}\label{ag} a^V = \sum_{\gamma \in \Gamma} \, a_{e^*_{\gamma}} \tp e_{\gamma} +
 \sum_{\gamma \in \Gamma}\,
 a_{e_{\gamma}} \tp e_{\gamma}^* =a^V_+ + a^V_- ,
\end{equation} where $a_{e^*_\gamma}, \, a_{e_\gamma} \in {\rm End}\, W.$ The
coefficients in (\ref{ag}) are defined via the trace (\ref{res}).
For any element $a^V \in {\rm End}\, W \tp V^{\g}$ and any $f \in
V^{\g}$ we can define $a^V_f \in {\rm End}\, W$ as \begin{equation} a^V_f =
Res_{V^\g}((1 \tp f )\cdot a^V ).\end{equation}

We have a natural action of $\CC[\d_{1},\ldots ,\d_{p}]$ on ${\rm
End}\, W \tp V^{\g},$ induced by the action on $V^{\g},$ so that
\begin{equation} (\d_i \, a)^V _f= - a^V _{\d_i f}, \end{equation} or, more general:
\begin{equation}
(\d_1^{k_1}\,\d_2^{k_2}\ldots \d_p^{k_p} \, a)^V _f= (-1)^{k_1 +k_2+
\ldots k_p} \, a^V _{
\d_p^{k_p}\,\d_{p-1}^{k_{p-1}}\ldots \d_1^{k_1}\,f}. \end{equation}

We will say that two formal distributions $a^V$ and $b^V$ are
mutually local if they commutator can be expressed as a finite
linear combination of the $\delta_V$ and its derivatives:
\begin{equation}\label{cr} [a^V_z, \, b^V_w]= \sum_{k_1,\ldots k_p } \,
c^V_{(\, \d_1^{k_1}\ldots \d_p^{k_p}v_{\mu})^* ;\,w}
\,\d_{1,w}^{k_1}\,\d_{2,w}^{k_2}\ldots \d_{p,w}^{k_p}\, \delta_V .
\end{equation}
 Since $ \delta_V \in V_z^{\g} \tp (V_w^{\g})^*$ in the previous
 formula, we need to specify to which factor we apply the operator
 $\d_k .$ $\, \d_{k,w}$ means that we apply it to the second
factor in the tensor product.

 The coefficient $c^V_{(\, \d_1^{k_1}\ldots \d_p^{k_p}v_{\mu})^*;
 \,w} \in {\rm End}\, W \tp V^{\g}_w$
in this decomposition can be viewed as the $v_{k_1,\ldots k_p
}$-product of fields $a^V$ and $\, b^V.$ The element
$v_{k_1,\ldots k_p }$ is in the space $V^*_\mu ,$ dual to some
representation of $\g$ with a lowest weight $\mu$ and the lowest
vector $v_\mu \in V_\mu ,$ and is dual to
$\d_1^{k_1}\,\d_2^{k_2}\ldots \d_p^{k_p}$ in the sense that
\begin{equation}
 \d_p^{k_p}\,\d_{p-1}^{k_{p-1}}\ldots \d_1^{k_1}
\,v_{s_1,\ldots s_p } =(-1)^{k_1 +k_2+ \ldots k_p}\, v_\mu^*
.\end{equation} In particular, $\d_i v_\mu^* =0.$

For any $f$ and $g \in V^\g$ we can define $[a^V_f ,\, b^V_g]\in
{\rm End}\, W \tp V_z^{\g} \tp V_w^{\g}$ as: \begin{equation} [a^V_f ,\,
b^V_g]=Rez_{w}Rez_{z}\,( [a^V, \, b^V]\cdot (1\tp f \tp1) \cdot
(1\tp1 \tp g)). \end{equation}

 To consider the space of formal distributions
with the values in a Lie algebra $\a$ we need to impose the
 skewsymmetry and the Jacobi identity axioms to the $a^V_{(v)}
b^V$-products of two fields, $\,\,v\in V^*_\mu .$ Denote by $R^\g$
the set of all local fields (distributions) on $V^\g$ with the
values in $\a.$
 The $v $-products define a $\CC$-linear map \begin{equation}\label{pr} \ca :
\,\, V^*_\mu \tp R^\g \tp R^\g \longrightarrow R^\g .\end{equation} This map
satisfies the following axioms, that are analogues of axioms
(C1)-(C2) in the Def.1.1:

(H1) For any $a, b \in R^\g$ the map $\ca : \,\, V^*_\mu \tp a \tp
b \longrightarrow R^\g$ is non zero only on a finite dimensional
subspace $V_{ab}^* \subset V^*_\mu ,$

(H2) $\ca \circ( 1 \tp \d_i \tp 1) = \ca \circ ( \d_i \tp 1\tp
1)\,\,$ and $\,\, \ca \circ (1\tp 1 \tp \d_i)= \d_i \circ \ca +
\ca \circ ( 1 \tp \d_i \tp 1).$
 \begin{definition}
 A(Lie) generalized conformal algebra
  $R^\g , $
 associated with the space $V^*_\lambda$ is a
 $\CC[\d_{1},\ldots ,\d_{p}]$-module, endowed with the map (\ref{pr}),
 satisfying axioms (H1)-(H2) as well as the skewsymmetry and the Jacobi
 identity axioms.
 A conformal algebra $R^\g$ is called $finite$ if $R^\g$ is a finitely generated
 $\CC[\d_{1},\ldots ,\d_{p}]$-module.
The $rank$ of conformal algebra $R^\g$
 is its rank as a
 $\CC[\d_{1},\ldots ,\d_{p}]$-module.
\end{definition}

Since $\n_+$ is a non commutative algebra, it is difficult to
write explicitly in terms of $\ca$ the skewsymmetry and the Jacobi
identity axioms. As for Conformal algebras these axioms have more
simple form in terms of the so-called $\lambda-bracket$ (see for
reference \cite{K2}) defined as: \begin{equation} [a_\la
b]=\sum_{n \in \ZZ}\la^{(n)} \, a_{(n)}b .\end{equation}

We define a non-commutative analogue of the $\lambda-bracket $ as
follow. Consider a ${\rm Poincar\acute e}$-Birkhoff-Witt basis
${\bf B}=\{e_{\al_1}^{k_1}\ldots e_{\al_p}^{k_p}\} \in U(\n_+).
\,$ Let ${\bar \lambda} = \la_1^{k_1}\ldots \la_p^{k_p}$ be the
symbol of $e_{\al_1}^{k_1}\ldots e_{\al_p}^{k_p}$ in
$\CC[\la_1,\ldots ,\la_p]$ and ${\bar \d }=\d_1^{k_1}\ldots
\,\d_p^{k_p}\,\,$ its image in $\, \CC[\d_1,\ldots ,\d_p].\,$
Denote by $\phantom{\,}^t {\bar \la} =(-\la_p)^{k_p}\ldots
(-\la_1)^{k_1}$ the image of ${\bar \la}$ under the transposition
map. Define the ${\bar \la}$-bracket of two elements $a,\, b \in
R^\g$ as: \begin{equation}[a_{\, \bar \la}\,
b]=\sum_{\d_1^{k_1}\ldots \d_p^{k_p} \in {\bf
B}}\,\la_1^{(k_1)}\ldots \la_p^{(k_p)} \, a_{(\, \d_1^{k_1}\ldots
\d_p^{k_p}v_{\mu})^* }\,b .\end{equation}

The axioms (H1)-(H2) are rephrased as follows:

(H1)$[a_{\, \bar \la}\, b] \in \CC[\la_1,\ldots ,\la_p] \tp R^\g,$

(H2) $[ \d_i a_{\, \bar \la}\, b]= -\la_i \,[a_{\, \bar \la}\,
b],\;\;\;\;\;[a_{\, \bar \la}\, \d_i b] =(\d_i +\la_i)\,[a_{\,
\bar
\la}\, b]$,\\
 the the skewsymmetry axiom:

(H3) $[a_{\, \bar \la}\, b]= -[b_{\, \phantom{i}^t(\bar \d + \bar
\la)}\, a],$\\
the Jacobi identity axiom:

(H4) $[a_{\, \bar \la}\,[b_{\, \bar \mu}\, c]= [[a_{\, \bar
\la}\,b]_{\,\bar \la+ \bar \mu}\, c]+[b_{\, \bar \mu}\,[a_{\, \bar
\la}\, c]].$

In this setting the axioms (C1)-(C4) of conformal algebra $R$ in
Def.1.1, applying to the affine Kac-Mody and Virasoro cases say
that there are only two opposite cases:

(i) $R/\d R$ is a Lie algebra and $V_\mu$ is one dimensional
module. Then for $a,b \in R/\d R$ we define $a_{(v^*_\mu)}b=[a,
\,b].$ This is the Current conformal algebras case.

(ii) $R/\d R$ is one dimensional (rank one conformal algebra) and
$V_\lambda$ is a self-dual $\CC[\d]$-module that has the structure
of Lie algebra.

More precisely in the second case we take $\mu = -\al \; $ where
$\al$ is the root of $\slt.$ Then, $V_\al \simeq \slt$ with $v_\al
=Y$ and the action $\d$ given by $\d=ad \,(X).$ Here $X,\,H,\,Y$
is the standard basis of $\slt .$ Normalize this basis as $v_{\al}
=Y, \, H= \d Y, \,-2X=\d^2 Y $ and $\d X= 0.$
 The dual basis is $X, \,\, H/2,\,\, -Y/2 .$
Take $L=Y$ and define the products, corresponding to the elements
of the dual basis as:\begin{equation}\label{vcon} L_{(0)}L =
L_{(X)}L=[X,\,L]=\d L, \,\;\; L_{(1)}L = L_{(H/2)}L=[-H,\,L] =
2L.\end{equation} The axiom (H2) reads as $ (\d L)_{(e)}L=L_{(\d e)}L ,$ in
particular $(\d L)_{(H/2)}L=L_{(\d H/2)}L= -L_{(X)}L .$ These
relations between the commutation relations of the $\slt$ Lie
algebra and n-products of the Virasoro conformal algebra is a one
more manifestation of the link between Vir and $\slt$ Lie
algebras.

 Given a generalized conformal algebra $R^\g , $
 a Lie algebra of $V^\g$-local formal distributions associated to it
 is defined as follows.
 Fix a basis $\{e_{\gamma}, \, e_{\gamma}^*\},\, \, \, \gamma \in \Gamma$
 in $V^\g .$ Consider a vector space over $\CC$ with the basis

 $a^V_{e_{\gamma}}, \, a^V_{e^*_{\gamma}},$ where $a^V \in R^\g ,$
Then the Lie algebra of $V^\g$-local formal distributions is a
quotient of this space by the $\CC$-span of all elements of the
form \begin{equation} (\lambda a^V + \mu b^V
)_{e_{\gamma}}-(\lambda a^V)_{e_{\gamma}} - \mu (b^V
)_{e_{\gamma}},\;\;\;\;\;\; (\lambda a^V + \mu b^V
)_{e^*_{\gamma}}-(\lambda a^V)_{e^*_{\gamma}} - \mu (b^V
)_{e^*_{\gamma}}; \end{equation} \begin{equation} (\d_i \, a)^V_
{e_{\gamma}}+ a^V _{\d_i {e_{\gamma}}},\;\;\;\;\;\;\;(\d_i \,
a)^V_ {e^*_{\gamma}}+ a^V _{\d_i {e^*_{\gamma}}},
\;\;\;\;\lambda,\; \mu \in \CC .\end{equation}

\section{ The space $ V^{\g}$ }

In this section we discuss the basic ideas about the realization
 of the space $ V^{\g}$ in general case for any
simple Lie algebra $\g.$ We identify the space $ V^{\g}$ with the
space of regular function on an open subset (manifold) of the flag
manifold $B_-\backslash G$, where $G$ is the simply-connected Lie
group corresponding to the Lie algebra $\g$.

As a vector space, the Lie algebra $\g$
  has a triangular decomposition
\begin{equation}
 \g =\n_{+} \oplus \h \oplus \n_-,
\end{equation} where $\h$ is the Cartan subalgebra and $\n_{\pm}$ are
the upper and lower nilpotent subalgebras. Let \begin{equation}
\b_{\pm}= \h \oplus \n_{\pm} \end{equation} be the upper and lower
Borel subalgebras. Let $N_{\pm}$ (respectively, $B_{\pm}$ be the
upper and lower unipotent subgroups (respectively, Borel
subgroups) of $G$ corresponding to $\n_{\pm}$ (respectively,
$\b_{\pm}$). Let $\{\al_i\}  \, \, i=1,..., l$ be the root basis
of $\g$ and the $\{h_i \}_{i=1,.., l},\;\;$ the coroot basis of
$\h$. We choose a root basis of $\n_{\pm}$ , $\{e^{\al}\}_{\al \in
\Delta_{\pm}}$, where $\Delta_{+}$ ($\Delta_{-}$) is the set of
positive ( negative) roots of $\g$, so that $[ h, \, e^{\al}] =
\al (h) e^\al$ for all $h \in \h.$

Consider the flag manifold $B_-\backslash G.$ It has a unique open
$N_+$-orbit, the so-called big cell $U=[1]\cdot N_+ \subset
B_-\backslash G$, isomorphic to $N_+$. Since $N_+$ is a unipotent
Lie group, the exponential map $\n_+ \to N_+$ is an isomorphism
 and we have
$U \simeq \CC^3 .$
 From the action of $N_+$ on $U$ we can introduce
a system $\{y_\al \}_{\al \in \Delta_+}$ of $homogeneous$
coordinates on $U$. The homogeneous means that
\begin{equation}\label{basis1} h \cdot y_\al = -\al (h) y_\al
\end{equation} for all $h \in \h$.

The action of $G$ on $B_-\backslash G$ gives us a map from $\g$ to
the Lie algebra of vector fields on $B_-\backslash G ,$ and hence
on its open subset $U\simeq N_+. $ Thus we obtain a Lie algebra
homomorphism $\g \to \rm{Vect}\, N_+$ With respect to this action
the space Fun $N_+$ of regular functions on $U$ has structure of
the contragradient Verma module $M^*_0$ with lowest weight $0$
(for more details se \cite{F1}). We have a natural pairing
$U(\n_+) \times \rm{Fun }\,N_+ \to \CC,$ which maps $(P\, , \, A)$
to the value of the function $P \cdot A$ at the identity element
of $N_+$, for any $A \in$ Fun $N_+$ and $P \in U(\n_+).$

The vector $1 \in \rm{Fun }\,N_+ $ is annihilated by $\n_+$ and
has weight $0$ with respect to $\h$. Hence there is a non-zero
homomorphism Fun $\,N_+ \to M^*_0$ sending $1 \in \rm{Fun }\,N_+$
to a non-zero vector $v^*_0 \in M^*_{0}$ of weight $0.$ Since both
$\rm{Fun} \,N_+$ and $ M^*_{0}$ are isomorphic to $U(\n_+)^{\vee}$
as $\n_+$-modules, this homomorphism is an isomorphism.

 It is known (\cite{F1}) that we can identify the Module $M^*_{\chi}$ with an arbitrary weight
$\chi$ with $\rm{Fun}\, N_+$, where the latter is equipped with a
modified action of $\g$. Recall that we have a canonical lifting
of $\g$ to $\cd_{\le 1} (N_+)$ of differential operators on $U$ of
order one, in a way that $a \to \xi_a.$ The modified action is
obtained by adding to
 each $ \xi_a$ a function
 $\phi_a \in \rm{Fun} \,N_+.$ The modified differential operators $\xi_a +\phi_a$
 satisfy the commutation
relations of $\g$ if and only if the linear map $\g \to \rm{Fun}\,
N_+.$ given by $a \to \phi_a$ is a one-cocycle of $\g$ with
coefficients in ${\rm Fun}\, N_+.$

If we impose the extra condition that the modified action of $\h$
on $\rm{Fun}\, N_+$ remains diagonalizable, we get that our
cocycle should be $\h$-invariant: $\phi _{[h, a]} = \xi_h \cdot
\phi_a,$ for all
 $\, \; h \in \h, \,\, a\in \g.$ The space of $\h$-invariant one-cocycle of $\g$
 with coefficients in $ \rm{Fun}\, N_+$
is canonically isomorphic to the first cohomology of $\g$ with
coefficients in $ \rm{Fun}\, N_+$ (see \cite{F1}). By Shapiro's
lemma we have $ H^1 (\g, \, \rm{ Fun} \,N_+.) = H^1(\g, \,M^*_0)
\simeq H^1 (\b_-, \, \CC_0)= (\b_-/[\b_-, \b_-])^* \simeq\h^*.$
Thus, for each $\chi \in \h$ we have an embedding $\rho_{\chi} \,
: \,\g \hookrightarrow \cd_{\le 1} (N_+)$
 and the structure of
$\h^*$-graded $\g$-module on $\rm{Fun}\, N_+.$

More detailed analyze shows that the action of $\n_+$ on $\rm{Fun}
\,N_+$ is not modified and the action of $h \in \h$ is modified by
$h \to h + \chi(h).$ This comes from the reason that the weight of
any monomial in $\rm{Fun}\, N_+ $ is equal to the sum of negative
roots (\ref{bas}) and the $\h$-invariance of the one-cocycle $
\xi_h \cdot \phi_{e^\al} =\phi _{[h, e^\al]} = \al (h)
\phi_{e^\al}, \, \, \, \al \in \Delta_+.$
 Therefore $ \phi_{e^\al}=0$ for all $\al \in\Delta_+$.
 The vector $1 \in \rm{Fun}\, N_+$ is still
annihilated ny $\n_+ ,$ but now it has weight $\chi$ with respect
to $\h.$ Hence there is a non-zero homomorphism $\rm{Fun }\,N_+
\to M^*_\chi$ sending $1 \in \rm{Fun}\,N_+$ to a non-zero vector
$v^*_\chi \in M^*_{\chi}$ of weight $\chi.$ Since both $\rm{Fun}
\,N_+$ and $ M^*_{\chi}$ are isomorphic to $U(\n_+)^{\vee}$ as
$\n_+$-modules, this homomorphism is an isomorphism.

For the reasons explained below we will consider the modified
action of $\g$ on $\rm{Fun} \,N_+$ with highest weight
\begin{equation} \chi = - \rho, \, \, \, \rm{ where }\, \, \, \,
\rho=\frac{1}{2}\sum_{\al\in \Delta_+} \al. \end{equation}

\begin{proposition}
With respect to the modified action of the Lie algebra $\g$
 the space $\rm{Fun} \,N_+$ has the structure of
the Verma module $M_{-\rho}$ with highest weight $-\rho.$
\end{proposition}

The space $\rm{Fun} \,N_+$ is closed under the multiplication and
$\n_+$ acts on this space by derivations. We will construct the
bigger space with the same action of $\g$ that has all the
properties of Definition 3.1.

 We have the natural action of the Weyl group
$W$ of $\g$ on the flag manifold $B_-\backslash G . $ From the
definition $W \simeq N(T^l)/T^l$, where $T^l$ is maximum torus in
$G$ and $N(T^l)$ is its normalizer.  We have a particular element
$w_{0} \in W$ called the longest element (see [\cite{W}]),
satisfying the following three conditions:

$(i) \, \, \,w_0(\Delta_+) = \Delta_- ,$

$(ii)\, \, w_0 \cdot \rho= - \rho ,$

$(iii) \, l(w_0)= |\Delta_+|,$\\
where $ \rho=\frac{1}{2}\, \sum_{\al\in \Delta_+}\, \al. \,$ This
element is uniquely defined, and satisfies $w_0^2=e.$ For example,
for type $A_l$ and any simple positive root $\al_i$, $\, w_0\al_i
= - \al_{l+1-i}, \, \, 1\le i \le l.$ It means that $w_0 N_+
w_0^{-1} =N_-.$

 Consider the big cell $U^{*} \in B_-\backslash G$ dual to $U= [1]
 \cdot N_+ \subset B_-\backslash G.$ We have
 \begin{equation}
U^* =[1] w_0 \cdot N_-  = [1] w_0 w_0 N_+ w_0 =  U w_0
.\end{equation} Let $\rm {Fun } \, U^*$ be the ring of regular
function on $U^*.$ We have $U^*\simeq N_-$.
 From the action of $N_-$ on $U^*$ we can introduce
a system $\{x_\al \}_{\al \in \Delta_+}$ of $homogeneous$
coordinates on $U^*$. Homogeneous means that
\begin{equation}\label{bas} h \cdot x_\al = \al (h) x_\al
\end{equation} for all $h \in \h$.

On the intersection $U\cap U^*$ we can consider
 the change of variables
on $U\cap U^*$ from $\{y_\al \}$ to $\{x_\al \}.$

 There is an ideal in ${\rm Fun }\, U^*$ that has
a structure of the Verma module $M^{\vee}_{\rho}$ with the lowest
weight $\rho$ with respect to the modified action of the Lie
algebra $\g .$

\begin{proposition}
There is an element $f \in \rm{Fun }\, U$ that has the following properties:\\
$\, \, \,(i) \, \, \, f $ has a weight $ \, -2\rho \,$ with respect to the action of $\h,$\\
$\,\, (ii) \,\,\phi= (f)^{-1}$ is an element of $\rm{Fun } \,U^*$ of weight $\rho,$\\
$\,(iii) \, \phi $ is annihilated by $\n_- ,$\\
$\,(iv) \, $ the ideal $\, V_\rho =\phi \cdot \rm{Fun } U^* \in
{\rm Fun } U^*$ has the structure of the Verma module $M'_\rho$
with the lowest weight $\rho$ with respect to the modified action
of the Lie algebra $\g.$

\end{proposition}

 We have a natural pairing $<\cdot,\, \cdot>$  between the Verma module $M_{-\rho}
\simeq \rm{Fun } \,U$ with highest weight $-\rho$ and the Verma
module $M'_\rho \simeq V_ \rho \in \rm{Fun } \,U^*$ with lowest
weight $\rho.$

Consider the space $V_0 = \rm{Fun }\, U\, \oplus \, \, \phi \cdot
\rm{Fun }\, U^*$. The element $\phi \in V_0$ define an invariant
trace on $V_0$ in a way:
\begin{equation}
Res \,f =\,<f,\,1>,\;\;\;\;f \in V_0
 \end{equation}

 This space have almost all properties of the
space from the Definition 2.1, with exception that it is not
closed under the multiplication. Now we will construct the bigger
space $V^{\g}$, that contains $V_0 $ as the subspace and closed
under the multiplication.

Conceder a complex submanifold $M \in B_-\backslash G$ that is the
intersection of $U$ and $U^* $
 \begin{equation} \label{m}M \,= U \, \cap \, U^* .
 \end{equation}
As a set $M$ is isomorphic to the complex space $\CC^n \setminus
D_\phi \; \;\;n=|\Delta_+|,$ where $D_\phi$ is the divisor defined
by $ \phi=\phi_1(y_{\al_1}, ..., y_{\al_n})\cdot ... \cdot
\phi_1(y_{\al_1}, ..., y_{\al_n})=0, \;\;i=1,...,l $ and each
divisor $D_{\phi_i} \,: \phi_i(y_{\al_1}, ..., y_{\al_n})=0$ is a
simple divisor.

 Denote by $V^{\g}$ the space $\rm{Fun} \, M$ of regular
functions on $M.$

\begin{theorem}
  The space  $V^{\g}$ has all the properties listed in the Definition
  2.1.
\end{theorem}

The proof to this theorem  the Proposition 4.2  and the explicit
construction for the manifold $M \,= U \, \cap \, U^* $ for the
Lie algebra $\g$ of one of the type $A_n,\,B_n,\,C_n, D_n$ will be
given in the future publications of the author. The case $\g=A_2$
is discussed in the next section.

\section{Explicit realization of the space $ V^{\slf}$ }

In this section we give an explicit construction of the space $
V^{\g}$ for $ \g=\slf $ as the space of regular function on an
open subset (manifold) of the flag manifold $B_-\backslash G$,
where $G=SL(2, \CC)$ is the simply- connected Lie group
corresponding to $\g=\slf$.
 For this section we assume that $ \g=\slf$.
 Fix a triangular decomposition
\begin{equation}
 \g =\n_{+} \oplus \h \oplus \n_-,
\end{equation} where $\h$ is a Cartan subalgebra and $\n_{\pm}$ are
the upper and lower nilpotent subalgebras. Let \begin{equation}
\b_{\pm}= \h \oplus \n_{\pm} \end{equation} be the upper and lower
Borel subalgebras. Let $N_{\pm}$ (respectively, $B_{\pm}$ be the
upper and lower unipotent subgroups (respectively, Borel
subgroups) of $G$ corresponding to $\n_{\pm}$ (respectively,
$\b_{\pm}$). Let $h_i = \al_i^{\vee}, \, \, i=1, 2$ be the $i$-th
coroot of $\g$. The set $\{h_i \}_{i=1, 2}$ is a basis of $\h$. We
choose the root basis of $\n_{\pm}$ , $\{e^{\al}\}_{\al \in
\Delta_{\pm}}$, where $\Delta_{+}$ ($\Delta_{-}$) is the set of
positive ( negative) roots of $\slf$, so that $[ h, \, e^{\al}] =
\al (h) e^\al$ for all $h \in \h.$

Consider the flag manifold $B_-\backslash G.$ It has a unique open
$N_+$-orbit, the so-called big cell $U= N_+ \cdot[1] \subset
B_-\backslash G$, isomorphic to $N_+ ,$ it means that  $U \simeq
\CC^3 .$
 From the action of $N_+$ on $U$ we can introduce
a system $\{y_\al \}_{\al \in \Delta_+}$ of $homogeneous$
coordinates on $U$ in a following way. Let
\[
n=\left(
\begin{array}{ccc}
1 & a & c\\
0 & 1 & b\\
0 & 0 & 1\\
\end{array}
\right)
\]
be an element of $N_+$ and $a, b, c \in \CC,$ then $(a, b, c)$ are
 homogeneous coordinates on $U$.

  We have a Lie algebra homomorphism
$\slf \to \rm{Vect}\, N_+$ With respect to this action the space
Fun $N_+$ of regular functions on $U$ has structure of the
contragradient Verma module $M^*_0$ with lowest weight $0.$
 We have a natural pairing $U(\n_+) \times
\rm{Fun }\,N_+ \to \CC,$ which maps $(P\, , \, A)$ to the value of
the function at the identity element of $N_+$, for any $A \in$ Fun
$N_+$ and $P \in U(\n_+).$ This pairing define a isomorphism
between the  $\rm{Fun} \,N_+$ and $ M^*_{0}.$

Consider the modified action of $\slf$ on $\rm{Fun} \,N_+$
described in the previous section  with the highest weight
\begin{equation} \chi = - \rho= -(\al_1+\al_2)=-\al_3, \end{equation}
where  $\{\al_1, \, \al_2,\,\, \al_3=\al_1 + \al_2 \}$ are
positive  roots for $\g =\slf .$ The corresponding basis in $\slf$
is $\{e_1, \, e_2,\, e_3, \, h_1,\, h_2,\, f_1,\, f_2,\, f_3\},$.
The explicit formulas for the modified $\slf$-action in terms of
homogeneous coordinates $a, b, c$ on $U$:
$$\xi_{\al_1}=\xi_{-\rho}(e_1) = \frac{\d}{\d a}=\d_1,$$
$$\xi_{\al_2}=\xi_{-\rho}(e_2) = \frac{\d}{\d b} + a \frac{\d}{\d c} = \d_2
,$$ $$ \xi_{\al_3}=\xi_{-\rho}(e_3) = \frac{\d}{\d c} =\d_3 ,$$
$$\xi_{h_1}=\xi_{-\rho}(h_1) = -2a \frac{\d}{\d a}+b\frac{\d}{\d b} -c \frac{\d}{\d c}-1
,$$
$$\xi_{h_2}=\xi_{-\rho}(h_2) = a \frac{\d}{\d a}-2b\frac{\d}{\d b} -c \frac{\d}{\d c}-1
,$$
$$\xi_{-\al_1}=\xi_{-\rho}(f_1) = -a^2 \frac{\d}{\d a}- (c-ab)\frac{\d}{\d b} -ac \frac{\d}{\d c}-a
,$$
$$\xi_{-\al_2}=\xi_{-\rho}(f_2) = c\frac{\d}{\d a} -b^2
\frac{\d}{\d b} -b ,$$ \begin{equation}\label{action}
\xi_{-\al_3}=\xi_{-\rho}(f_3) = -ac\frac{\d}{\d a}- b(c-ab)\frac{\d}{\d b} -c^2 \frac{\d}{\d c}-(2c-ab).\\
\end{equation}
\begin{proposition}
With respect to the action of the Lie algebra $\slf$ given by
(\ref{action}), the space $\rm{Fun} \,N_+$ has the structure of
the Verma module $M_{-\rho}$ with highest weight $-\rho.$
\end{proposition}

The space $\rm{Fun} \,N_+$ is closed under the multiplication and
$\n_+$ acts on this space by derivations. We will construct the
bigger space with the same action of $\slf$ that has all the
properties of Definition 3.1.

 We have the natural action of the Weyl group
$W$ of $\slf$ on the flag manifold $B_-\backslash G . $ From the
definition $W \simeq N(T^l)/T^l$, where $T^l$ is maximum torus in
$G =SL(3, \CC)$ and $N(T^l)$ is its normalizer. In the case of
$\slf$ the Weyl group $W \simeq S_3$, where $S_3$ is a group of
permutations.
 The longest element  $w_{0} \in W$ corresponds to the permutation
 $(1,\,3)$ via this identification.

 Consider the big cell $U^{*} \in B_-\backslash G$ dual to $U=
[1]\cdot N_+  \subset B_-\backslash G.$ We have
 \begin{equation}
U^* = [1] w_0 \cdot  N_-= U w_0  .\end{equation} Let $\rm {Fun }
\, U^*$ be the ring of regular function on $U^*.$ We have
$U^*\simeq N_-$.
 From the action of $N_-$ on $U^*$ we can introduce
a system $\{y_\al \}_{\al \in \Delta_-}$ of $homogeneous$
coordinates on $U^*$. Let
\[
n=\left(
\begin{array}{ccc}
1 & 0 & 0\\
s & 1 & 0\\
t & r & 1\\
\end{array}
\right)
\]
be an element of $N_-$ and $s, r, t \in \CC .$ Then $s, r, t $ are
a homogeneous coordinates on $U^* .$  With respect to the action
of the Cartan subalgebra (\ref{action}) any monomial of the form
$s^n r^m t^k $ has the weight equal to the $-\rho+n\al_1 + m \al_2
+k \al_3.$

The intersection $U\cap U^*$ is an open subset in $B_-\backslash
G.$ On the $U\cap U^*$ the relations between the homogeneous
coordinates $a, b, c$ and $s, r, t $ are given by:
\begin{equation}\label{zamena} a= \frac{r}{t},\, \, \, \,b= - \frac{s}{t-rs},
\,\, \, \, c =\frac{1}{t}; \,\, \,\, \rm{or}\, \, \,\,\,
s=-\frac{b}{c-ab}, \,\,\,\, r=\frac{a}{c}, \,\,\,\, t=\frac{1}{c}.
\end{equation}

We will show that there is an ideal in ${\rm Fun }\, U^*$ that has
a structure of the Verma module $M^{\vee}_{\rho}$ with the lowest
weight $\rho$ with respect to the action of the Lie algebra $\g$
given by (\ref{action}).

\begin{proposition}
There is an element $f(a,b,c) \in \rm{Fun }\, U$ that has the following properties:\\
$\, \, \,(i) \, \, \, f $ has a weight $ \, -2\rho \,$ with respect to $\h,$\\
$\,\, (ii) \,\,\phi= (f)^{-1}$ is an element of $\rm{Fun } \,U^*$ of weight $\rho,$\\
$\,(iii) \, \phi $ is annihilated by $\n_- ,$\\
$\,(iv) \, $ the ideal $\, V_\rho =\phi \cdot \rm{Fun } U^* \in
{\rm Fun } U^*$ has the structure of the Verma module $M'_\rho$
with the lowest weight $\rho$ with respect to the action of the
Lie
algebra $\g$ given by (\ref {action}),\\
$(v)$ the differential form $\frac{da\wedge db \wedge
dc}{c(c-ab)}$ is invariant under the transformation
(\ref{zamena}): $ \frac{da\wedge db \wedge
dc}{c(c-ab)}\longrightarrow \frac{ dr \wedge ds \wedge
dt}{t(t-sr)}$ .
\end{proposition}

{\bf{Proof.}} Consider a function $f(a, b, c) = c(c-ab) =
({t\,(t-rs)})^{-1}.$ Straightforward calculations show that this
element have all the properties, listed above.

We have a natural pairing  $<\cdot ,\, \cdot>$  between the Verma
module $M_{-\rho} \simeq \rm{Fun } \,U$ with highest weight
$-\rho$ and the Verma module $M'_\rho \simeq V_ \rho \in \rm{Fun }
\,U^*$ with lowest weight $\rho.$

Consider the space $V_0 = \rm{Fun }\, U\, \oplus \, \, \phi \cdot
\rm{Fun }\, U^*$. This space have almost all properties of the
space from the Definition 2.1, with exception that it is not
closed under the multiplication. Now we will construct the bigger
space $V^{\g}$, that is the algebraic closer (with respect to the
multiplication) of the space $V_0 $ and we will show that
$V^{\g}\simeq {\rm Fun}\, U \, \cap \, U^*.$ This explains the
chose of the highest wight $\chi =-\rho$ for the modified action
of $\g$ on $\rm{Fun }\, N_+ .$

Conceder a complex submanifold $M \in B_-\backslash G$ that is the
intersection of $U$ and $U^*$
 \begin{equation} \label{m}M \,= U \, \cap \, U^* .
 \end{equation}
 Denote by $V^{\g}$ the space $\rm{Fun} \, M$ of regular
functions on $M.$ In order to obtain a nice description of $M$ and
of the space $V^{\g}$ we identify a point (a flag) in the flag
manifold $ B_-\backslash G$ with two orthogonal projective
vectors:
\begin{equation} \{(z_1: z_2 : z_3), \, \, \, (w_1: w_2 : w_3) \,|
\, \, \sum_{i =1}^{3} z_i w_i = 0\}.\end{equation}
 We have $[1] \in B_-\backslash G =\{(0, \,0,\, 1); \, (1,\,0,\,
0)\}$ and \begin{equation} U=\{(z_1: z_2 : z_3), \, \, \, (w_1: w_2 : w_3) | \,
\, z_3 \neq 0, w_1 \neq 0 \}. \end{equation} In terms of homogeneous
coordinates $(a,\, b,\,c )$ we have: \begin{equation}
 U= \{(-c+ab, \, -b,
\, 1), \, \, \, (1,\, a,\,c ), \,\,\, a, b, c \in \CC \}. \end{equation}

The natural action of the Weyl group $W $is just the corresponding
permutation of coordinates. If $s \in W \simeq S_3,$ then
\begin{equation} U_s =\{(z_1: z_2 : z_3), \, \, \, (w_1: w_2 :
w_3) | \, \, z_{s(3)} \neq 0, w_{s(1)} \neq 0 \} \end{equation}
where $U_s$ is the big cell in $B_-\backslash G ,$ that is the
image of $U$ under the action of the element $ s \in W$. In
particular, for the dual (opposite) sell we have $U^* = U_{w_0}$.
Thus,
\begin{equation}\label{proj}
 M=\{(z_1: z_2 : z_3), \, \, \, (w_1: w_2 : w_3) \,| \,
\, \sum_{i =1}^{3} z_i w_i = 0, \,\,{\rm{and}}\,\, z_i \neq 0,
\,\,w_i \neq 0, \,\, i=1,\,3 \} \end{equation} is the intersection of $\CC^*
\times \CC \times \CC^* \times \CC $ with the hypersurface
$\sum_{i =1}^{3} z_i w_i = 0$ in $\CP^2 \times \CP^2.$

If we fix homogeneous coordinates \begin{equation}
 (a,\, b,\,c ),\, \, \, a =\frac{w_2}{w_1}, \,\, c=
 \frac{w_3}{w_1}, \,\, b= -\frac{z_2}{z_3}
 \end{equation}
on $ U \simeq \CC^3,$ then $M$ is the total space of non trivial
bundle with the base $ (\CC)^2_{a,b}$ - the two-dimensional
complex space and the fiber $\CC^*_c\setminus\{c=ab\}:$

\begin{equation} M \stackrel{\CC^*_c \setminus\{c=ab\}}\longrightarrow
(\CC)^2_{a,b} . \end{equation}

The action (\ref{action}) of the $\n_+$ on $\rm{Fun}\, M$ in terms
of projective coordinates $(z_1: z_2 : z_3)$ and $ \, \, (w_1: w_2
: w_3)$ is: \begin{equation}\label{action1} \xi_{\al_1}= -z_2
\frac{\d}{\d z_1} + w_1 \frac{\d}{\d w_2},\;\;\;\; \xi_{\al_2}=
-z_3 \frac{\d}{\d z_2} + w_2 \frac{\d}{\d w_3},\;\;\;\;
\xi_{\al_3}=-z_3 \frac{\d}{\d z_1} + w_1 \frac{\d}{\d w_3},
\end{equation} and the element $\phi$ from the Proposition 5.2
  is
\begin{equation} \phi=\frac{z_3}{ z_1}\cdot
\frac{w_1}{w_3}.\end{equation} This element define
 $\n_+$-invariant trace on $V^{\g} \simeq \rm{Fun} \, M$
in the following sense. The space $V^{\g}$ is the space of all
regular function on $M$ or all functions on $\CC^3$ with the only
poles at $c=0$ and $c=ab$. We can choose a basis $\{e_{nml};\;
f_{nml}\}$ in $V^{\g}$ as:\begin{equation}\label{ba}
 e_{nml} = a^n c^m (c-a b)^l, \,\,\,n \ge
0,\;\;\;\,m,\,l \in \ZZ; \;\;\;\;\; f_{n ml} = b^n c^m (c-a b)^l,
\,\,\,\, n>0,\,\,m,l \in \ZZ
\end{equation}

\begin{proposition} Every regular function on $M$ is a unique linear combination
of the monomials of the basis (\ref{ba}).
\end{proposition}

Proof follows immediately  from the description of M in terms of
projective coordinates (\ref{proj}).

 The element
 \begin{equation}
 \phi=\frac{z_3}{ z_1}\cdot \frac{w_1}{w_3}= e_{0, -1, -1}
 \end{equation}
 is from this basis and for any regular function $f$
 on $M$ we can define a trace:\begin{equation}\label{tre}
 Res \, f =\, {\rm coefficient \;\;\;at} \,\, \phi \;\;\;{\rm in
 \;the\;\,
 decomposition\;\; with\;\; respect\;\; to \;\;this\;\; basis}.
 \end{equation}
 This trace is invariant with respect to the action of $\n_+$
 in a sense that
 \begin{equation}
   Res \,(\xi_{\al_i} (f))=0, \,\,\, {\rm for} \, \, \, i=1, 2, 3.
 \end{equation}

This invariant trace defines a bilinear invariant form on $V^\g
\simeq {\rm Fun} \, M$ as: \begin{equation}\label{form} < \, f, \,\, g \,
>= Res \, f\cdot g \;\;\;\;\;{\rm for}\;\;\;\; f,\, g \in {\rm Fun} \, M.\end{equation}

\begin{remark}
This invariant form, being restricted to the subspace $V_0 = Fun
\, U \oplus \phi \cdot Fun \, U* \subset V^{\g}$ is the natural
pairing between the Verma module $M_{-\rho} \simeq \rm{Fun } \,U$
with highest weight $-\rho$ and the Verma module $M'_\rho \simeq
V_ \rho \in \rm{Fun } \,U^*$ with lowest weight $\rho.$
\end{remark}

Consider the subspace $V_+$ in $\rm{Fun} \,M$ that consists of all
regular functions on $M_+ =(\CC)^{3}_{a,b,c} \setminus
D_{\{c=0\}}$. As a linear space $M_+$ is spanned by the elements
$\{e_{nml},\;n, \,l \ge 0,\,m\in \ZZ; \;\; \,f_{nml}, n>0,\;m
\ge,\;\;m\in \ZZ 0 \} $ of the basis (\ref{ba}). For simplicity,
we will denote the basis in $V_+$ by $\{e_{\gamma}\}, \,\,\gamma
\in \Gamma$, where $\Gamma$ is the set of all indexes of
$\{e_{nml},\, \{f_{nml}\} $. We assume that $e_0 = 1$. We have
\begin{equation} <e_{\gamma} \, , \phi \, > = 0 \;\;\;{\rm
for}\;\;\;\; \gamma \neq 0 \end{equation} and
\begin{equation} <e_{0} \, , \phi \, > = 1 .\end{equation}

 The above can be summarized in the following proposition.
\begin{proposition}(i) The space $V_+$ contains the unit function $1$ and it
is closed under the multiplication.\\
(ii)  $V_+$ is invariant under the action (\ref{action1}) of $\n_+
.$ \\
(iii) The bilibear form (\ref{form}) is non degenerated on $V^\g
.$\\
 (iiii)  $V_+$ is  the maximal isotropic subspace in $V^{\g}$
with respect to the bilinear form (\ref{form}).
\end{proposition}

Let $V_- \in \rm{Fun} \,M$ be the subspace dual to $V_+$ in
$V^{\g} $ with respect to the bilinear form (\ref{form}) and let
$\{e^*_{\gamma}\}, \,\,\gamma \in \Gamma$ be the basis in $V_-$
dual to $\{e_{\gamma}\}, \,\,\gamma \in \Gamma$. We have $e^*_{0}
= \phi$.

The space $V^{\slf}$ of all regular functions on $M$ is the direct
sum of two mutually dual subspaces $V^{\slf} = V_+ \oplus V_- .$
 The space $V^{\slf}$ satisfies to all properties (V1)-(V2) of the
Definition 2.1. So we can introduce the generalized formal delta
function associated with the space $V^{\g}$ as an element of the
space
\begin{equation}
  V^{\g} \tp (V^{\g})^* = V_+ \tp V_- \oplus V_- \tp V_+
\end{equation}
in a way described in the sec.2:
\begin{equation}
  \label{del1}
  \delta_{\scriptstyle{V^{\g}}} =\sum_{\gamma \in \Gamma} e_{\gamma} \tp
  e_{\gamma}^* \, \oplus \, e_{\gamma}^* \tp e_{\gamma}.
\end{equation}

 The space $V^{\g}$ is realized as the space of function
$\rm{Fun}\,M$ or, more general, as the space of formal
distributions on $M$ and ${\bf z} =(z_1, \, z_2, \, z_3)$ are
global coordinates on $M.$ Here $ z_1=a, \,\, z_2=b, \,\, z_3=c$,
where $(a,\,\,b, \, c)$ are homogeneous coordinates on $M.$
 The space $V^{\g} \tp (V^{\g})^*$ can be identified with the
space of distributions in two sets of coordinates ${\bf z} =(z_1,
z_2, z_3)$ and ${\bf w}=(w_1, w_2, w_3).$ In this case we will use
the notation \begin{equation} \delta_{V^{\g}}\,\, ({\bf z-w}), \end{equation} or for
short, $\delta_{V}\, ({\bf z-w}). $ We also can define
$\delta_{V}\, ({\bf z-w})_{\pm}$ as \begin{equation}\label{del+} \delta_{V}\,
({\bf z-w})_{+}=\sum_{e_{\gamma} \in V_+} e_{\gamma}^* \tp
e_{\gamma} \in V_-^z \tp V_+^w\end{equation} and \begin{equation} \delta_{V}\, ({\bf
z-w})_{-}=\sum_{e_{\gamma} \in V_+} e_{\gamma} \tp e_{\gamma}^*
\in V_+^z \tp V_-^w .\end{equation}

 So defined generalized
formal delta function has the most of the properties of the
standard delta function with respect to the trace (\ref{res}) and
differentiations from $\mathfrak{n}_+$: (a) For any formal
distribution $f({\bf z}) \in V^{\g}$ one has:
\begin{equation} Res_{ z} \,f({\bf z})\, \delta ({\bf z-w})=f({\bf w}), \end{equation}
 \begin{equation}
 (z_i -w_i)\,\delta ({\bf z-w}) =0 ;\end{equation}
(b) for any element $a \in \n_+$ one has: \begin{equation}
(\xi_a^z)^j \, \delta ({\bf z-w}) =(-
  \xi_a^{w})^j \, \delta ({\bf z-w}), \end{equation} where $\xi_{a_{z}}$ is the
image of $a$ in $\rm{Vect}\, M$.

 We would like to remark that for the Lie algebra $\slt$ the
 complex manifold $M \simeq U \bigcap U^* \simeq \CC^*$ is
 contactable to the one dimensional real manifold $S^1$.
 The trace (\ref{tree}) is the integral over the $S^1: \mid z \mid =1 $: \begin{equation}
 Res_z \; a(z) = \oint_{\mid z \mid =1} \; a(z)\, dz
\end{equation}
and the  formal Cauchy formula can be written in the form
\begin{equation}
 Res_z \; a(z) \partial_z^{(k)} \frac{1}{z-w}= (-1)^{k} \partial^{(k)} a(w)_+,
 \;\;\;\;{\rm for}\;\;\; \mid z \mid > \mid w \mid
\end{equation}
\begin{equation}
 Res_z \; a(z) \partial_z^{(k)} \frac{1}{z-w}= (-1)^{k+1}  \partial^{(k)} a(w)_-,
 \;\;\;\;{\rm for}\;\;\; \mid z \mid <  \mid w \mid
\end{equation}

Amazingly, we have very similar situation for the space $M^{\slf}
.$ The expression of the trace (\ref{tree}) in terms of the
integral over a real 3-dim. cycle as well as higher dimensional
version of the formal Cauchy formula will be given in
{\cite{G-K1}}
\begin{remark}
For the Lie algebra $\slt$ we have only two orbits of the Weyl
group in the Flag manifold. The space $\CC^*$ is the intersection
of two existing  dual orbit.  When we came to $\slf$ we have two
possibilities for generalization. One is to define the space
$M^\g$ as an intersection of two mutually dual orbits, what we we
did above, and the another one is to define $M^\g$ as
\begin{equation}
\tilde{M} \,= \,\bigcap_{w\in W} \, \, U_w . \end{equation} the
intersection of all orbits of the Weyl group in $ B_-\backslash
G.$
 Since the
action of $s \in W \simeq S_3$ is just the corresponding
permutation of coordinates: \begin{equation} U_s =\{(z_1: z_2 : z_3), \, \, \,
(w_1: w_2 : w_3) | \, \, z_{s(3)} \neq 0, w_{s(1)} \neq 0 \}. \end{equation}
Thus, \begin{equation} \tilde{M}=\{(z_1: z_2 : z_3), \, \, \, (w_1: w_2 : w_3)
\,| \, \, \sum_{i =1}^{3} z_i w_i = 0, \,\,{\rm{and}}\,\, z_i \neq
0, \,\,w_i \neq 0, \,\, i=1,\,2,\,3 \} \end{equation} is the intersection of
four-dimensional complex torus with the hypersurface \\$\sum_{i
=1}^{3} z_i w_i = 0$ in $\CP^2 \times \CP^2.$

In terms of homogeneous coordinates \begin{equation}
 (a,\, b,\,c ),\, \, \, a =\frac{w_2}{w_1}, \,\, c=
 \frac{w_3}{w_1}, \,\, b= -\frac{z_2}{z_3}
 \end{equation}
on $ U_1 \simeq \CC^3,$ $\tilde{M}$ can be described as the total
space of non trivial bundle with the base $ (\CC^*)^2_{a\,b},$ the
two-dimensional complex torus, and the fiber
$\CC^*_c\setminus\{c=ab\}:$ \begin{equation} \tilde{M}
\stackrel{\CC^*_c\setminus\{c=ab\}}\longrightarrow (\CC^*)^2_{a
\,b} . \end{equation} By the construction this space is invariant
under the action of the Weyl group $W.$

One of reason for the first choice is motivated by the following
theorems
\end{remark}
\begin{theorem}
The cohomology group \begin{equation} H^*_3(M^{\slf}, \;
\ZZ)\simeq \ZZ,
 \end{equation}
and
\begin{equation} H^*_3( \tilde{M}, \; \ZZ)\simeq
\ZZ\oplus \ZZ\oplus \ZZ .
 \end{equation}
\end{theorem}

\section{Generalized Affine Kac-Moody algebras associated with the space $
V^{\slf}$}

Let $\g$ be a simple Lie algebra over $\CC$. Consider the formal
loop algebra $L\g = \g((t))$. The affine algebra $\hg$ is defined
as the central extension of $L\g .$ As a vector space, $\hg =L\g
\oplus \CC K,$ and the commutation relations are: $[\, K, \, \cdot
] =0,$ and \begin{equation}\label{lie} [A \tp f(t), \, B \tp g(t)] = [A, \, B]
\tp f(t) g(t) + (A, \, B) Res_t (f \cdot g') \,K \end{equation} where
$(\cdot ,\, \cdot )$ is an invariant bilinear form on $\g ,$
normalized as in {\cite{K1}} by the requirement that $(\al_{max}
,\, \al_{max} )=2$.

Let $\{J^a\}_{a=1,\ldots dim \g}$ be a basis of $\g .$ Denote by
$J^a_n = J^a\tp t^n \in L\g,$ then $J^a_n, \, n \in \ZZ,$ and $K$
form a topological basis for $\hg .$ Consider a field (generating
function): \begin{equation} J^a (z)= \sum_{n \in \ZZ} \, J^a_n \tp
z^{-n-1}= J^a \tp \delta(t-z). \end{equation} For any element
$a\in \g$ we will define the corresponding field as
\begin{equation} a(z) = a \tp \delta(t-z).\end{equation} Here
$\delta(t-z)$ is the formal delta function associated with the
space ${\rm Fun}\, \CC^* \simeq V^{\slt}$. For any $a, \, b \in
\g$ the fields $a(z)$ and $b(z)$ are mutually local. The
commutation relations (\ref{lie}) in terms of the fields are:
\begin{equation}\label{loc} [a(z), \, b(w)]= [a,\,b] (w)\, \delta(z-w)+ K \cdot(a,
\, b) \, \d_w \delta(z-w). \end{equation}

In the same way we can define the affine Lie algebra, associated
with the space $V^{\slf} \simeq {\rm Fun} \, M.$ Let
\begin{equation} \delta_{\scriptstyle{V }}
=\delta_{\scriptstyle{V^{\slf}}}({\bf t-z}) \end{equation}
 be the delta function associated with the space ${V^{\slf}},$
 defined as (\ref{del1}).

 For any element $a \in \g$ we can
 define a field $a^{\slf} ({\bf z})$ on $M \,= U \, \bigcap \, U^*$ (\ref{m}) as:
 \begin{equation}
a^{\slf} (\bf{z})= a \tp \delta_{\scriptstyle{V }}\,(\bf{t}
-\bf{z}). \end{equation}

Consider the formal loop algebra $L^{\slf}\g =\g \tp {\rm Fun }\,
M$, associated with the manifold $M$ with the obvious commutation
relations.

 Let
 $ \xi_{i}^{z} = \xi_{\al_i}^{z}$ be
the image of $e_i, \, \, \, i=1,\,2,\,3$ in ${\rm Vect}\, M$
(\ref{action}). Then
 $ \d^z_i = \xi_{i}^{z}$ for each $i=1,\,2,\,3$ define a non trivial
 cocycle $c_i$ on $L^{\slf}\g$ as:
 \begin{equation}\label{coc}
c_i (A\tp f({\bf{t}}), \, B\tp g({\bf t})) =(A, \, B)\, Res_{{\bf
t}} (f \cdot \d_{i}^{t}( g)). \end{equation}
\begin{remark}
The Lie algebra $L^{\slf}\g$ has infinite dimensional central
extension. Let $f_i \in V^{\slf},\;\;i=1,\,2,\,3$ be  functions
such that $\sum_{i=1}^{3} \;\d_i f_i =0,$ then the vector field
$\sum _{i=1}^{3}\, f_i\d_i$ defines a non trivial central
extension of $L^{\slf}\g$ as:
 \begin{equation}\label{coc}
c_{f_1\,f_2,\,f_3}\, (A\tp f({\bf{t}}), \, B\tp g({\bf t})) =(A,
\, B)\, Res_{{\bf t}} (f \cdot \, \sum _{i=1}^{3}\,f_i \d_{i}^{t}(
g)).\end{equation} We can define the $\hgv$ as a infinite
dimensional central extension of the Lie algebra $L^{\slf}\g.$
\end{remark}

Not all these central extension are equally good from the point of
view of the conformal dimension.  The left hand side and the
right hand side of (\ref{loc}) have the same conformal dimensions.
Multiplication by $\d_w^k \delta(z-w)$ add $k+1$ to the conformal
dimension $\Delta$ of the field $a(z).$ The conformal dimension is
defined from the action of the Virasoro operators, in particular,
for the field  $a(z) $ of the conformal dimension $\Delta$ we have
 \begin{equation}[\,L_0,\; a(z)]= (z\,\frac{d}{dz}+\Delta)\,a(z).\end{equation}
  For the Affine
algebras all fields have conformal dimension $\Delta=1$ same as
the delta -function. For the Generalized Affine Kac-Moody algebras
associated with the space $ V^{\slf}$ we would like to consider
the central extensions compatible with the conformal dimension.

 Consider two grading operators on
$V^{\slf}$: \begin{equation}\label{grad} L_c = -c\, \d_3 - a \,
\d_1, \,\,\,\,\,\, L_{c-ab}= -(c-ab)\, \d_3 -b\, \d_2 .
\end{equation}
Here ${\bf z}=(a,\,b,\,c)$ where $(a,\,b,\,c) $ are homogeneous
coordinates on $U.$
\begin{proposition} Operators $L_c$ and $L_{c-ab}$ commute
and define  $\ZZ \times \ZZ$ grading  on $V^{\slf}$:
\begin{equation}\label{grading}
 V^{\slf}= \bigoplus_{(n_1,n_2)\in \ZZ}\;\;V_{n_1 n_2}\end{equation}
 where
$V_{n_1n_2}=\{f\in V^{\slf}| -L_c(f)=n_1 f;\;\;\,\,\,\, -L_{c-ab}
(f)= n_2f\}.$
\end{proposition}

Any element of $V^{\slf}$ is a linear combination of the
 elements of the basis (\ref{ba}), and every element from this
 basis have a form \begin{equation} a^n b^m c^k (c-ab)^l\end{equation}
 for some $n, m\in \ZZ_+$ and $ k,
 l \in \ZZ.$ From the explicit formulas for $\d_1, \, \d_2,\,\d_3 $ (\ref{action})
we have \begin{equation} a^n b^m c^k (c-ab)^l \in V_{n+k+l, \,
m+k+l}.\end{equation} This compleat the proof.

The grading operators $L_c $ and $ L_{c-ab}$ are analogues of the
energy operator $L_0$ that is a part of the Virasoro algebra
\begin{equation} Vir =\bigoplus_{n\in \ZZ}\,\, \CC L_n \oplus \CC
c .\end{equation} It is well known that that there is close
relations between the Virasoro and Affine Lie Algebras in many
aspects including the representation theory. In the next section
we  construct a natural higher dimensional analogue of the
Virasoro algebra associated with the space $V^{\slf}.$

We have the commutation relations:\begin{equation} [L_c,\;
a^{\slf} ({\bf z}]=(c\, \d_3 +a\,\d_1 +2)\;a^{\slf} ({\bf z});
\;\;\;\;\;\; [L_{c-ab},\; a^{\slf} ({\bf z})]=((c-ab) \,\d_3
+b\,\d_2 +2)\;a^{\slf} ({\bf z}).\end{equation} It means that the
field $a^{\slf}({\bf z}) $ has conformal dimension
$\Delta=(2,\,2)\in \ZZ\times \ZZ$ The generalized delta function
$\delta_{\scriptstyle{V }}\,(\bf{z} -\bf{w})$ as well has
conformal dimension $(2,\;2).$ From this we have the following
proposition:
\begin{proposition} Consider two cocycles on $L^{\slf}\g =\g \tp {\rm Fun
}\,M$ defined as:\begin{equation} c_i (A\tp f({\bf{t}}), \, B\tp
g({\bf t})) =(A, \, B)\, Res_{{\bf t}} \,(f \cdot L_{i}^{t}( g)),
\end{equation}
where ${\bf t}=(a,\, b,\, c)$ and $L_1=
\frac{1}{c(c-ab)}\;L_c,\;\;L_2=  \frac{1}{c(c-ab)}\;L_{c-ab}.$
Then the central extension defined by each of these cocycle has
the right conformal dimension with respect to the action of the
grading operators $L_c$ and $\;L_{c-ab}.$
\end{proposition}

\begin{remark}
 We have two Affine Kac-Moody subalgebras $\g \tp \CC[c,\,c^{-1}]$ and
 $\g \tp \CC[(c-ab),\,(c-ab)^{-1}]$ naturally imbedded in $L^{\slf}\g.$
 The restriction of each of these two cocycles to these
 subalgebras gives the standard Kac-Moody cocycle as in
 (\ref{loc}).
\end{remark}

\begin{definition} The Generalized Affine Kac-Moody algebra $\hgv$ associated with
the space $V^{\slf}$ is defined as a two-dimensional central
extension of $L^{\slf}\g.$ As a vector space, $\hgv =
L^{\slf}\g\oplus \CC K_1 \oplus \CC K_2 ,$ and the
 commutation relations: $[K_i,\, \cdot ]=0, \,i=1,\,2$ and
\begin{equation}\label{rel}
[a^{\slf} ({\bf z}), \, b^{\slf }({\bf w})]= [a,\, b]^{\slf} ({\bf
w})\, \delta_{\scriptstyle{V}}({\bf z -w}) +\sum_{i=1}^{2}\,K_i
\cdot(a, \, b) \,
 L_{i}^{w}\, \delta_{\scriptstyle{V}}(\bf z -w). \end{equation}
\end{definition}

 The commutations relations (\ref{rel}) give the description of
generalized Affine Kac-Moody algebras associated with the space
$V^{\slf}$ in terms of local fields. Here "local" means that the
commutator is annihilated by some polynomial function of $(z_i
-w_i).$

It is not completely appropriate to call the constructed algebras
the generalized Affine Kac-Moody algebras. These algebras have
many of the properties of "the generalized Affine Kac-Moody
algebras". The basic idea of R. Borcherds is to think of
generalized Affine Kac-Moody algebras as
 infinite dimensional Lie algebras which have most of the good
properties of finite dimensional
 reductive Lie algebras.

The original definition of generalized Affine Kac-Moody algebras,
given by R. Borcherds (\cite{B1}), is:

Consider a Lie algebra $G$ that has the following properties:

1. $G$ has an invariant symmetric bilinear form $(\,,\,)$.

2. $G$ has a (Cartan) involution $\omega$.

3. $G$ is graded as $G =\oplus_{n \in \ZZ} \, G_n$ with $G_n$
finite dimensional and with $\omega$ acting as $-1$ on the "Cartan
subalgebra" $G_0$.

4. $(g, \, \omega(g)) >0$ if $g \in G_n, \, g \neq 0$ and $n\neq
0.$

The Generalized Affine Kac-Moody algebras $\hgv$ associated with
the space $V^{\slf}$ are defined by the same conditions with one
small change: we replace condition 3 by

3* $G$ is graded as $G =\oplus_{{\bf n} \in \scriptstyle{\ZZ^k}}
\, G_{\bf n}$ with $G_{\bf n}$ not necessarily finite dimensional
and with $\omega$ acting as $-1$ on the "Cartan subalgebra" $G_0$.
Here $k=2$ is the rank of $\slf$.

In order to proof this statement we need to construct an invariant
symmetric bilinear form $(\,,\,)$ and
 a (Cartan) involution $\omega$ on $\hgv$.

The normalized invariant form $(\,,\,)$ on $\hgv$ can be described
as follows. Take the normalized invariant form $(\,,\,)$ on $\g$
and extend $(\,,\,)$ to the whole $\hgv$ by \begin{equation} (A\tp
f({\bf{t}}), \, B\tp g({\bf{t}}))= (A, \, B) (Res_t \,\, \phi
\cdot f \cdot g);\;\; (K_i, \,K_j)=0; \, \, \, \, \, (L^{\slf}\g,
\, \,\CC K_1 \oplus \CC K_2 \oplus \CC K_3)=0.\end{equation}

To define a Cartan involution of $\hgv$ consider the Cartan
involution $\tilde\omega$ of $\g$ and the transformation $w_0$ on
$M$ given by \begin{equation} a\to -\frac{b}{c-ab},\, \, \, \,b\to
 \frac{a}{c}, \,\, \, \, c \to\frac{1}{c}, \,\, \,\,
 \end{equation}
where ${\bf t}=(a,\,b,\,c)$ are homogeneous coordinates on $U .$
This transformation is defined by the change  from the homogeneous
coordinates on $U $ to the homogeneous coordinates on $U^*,$ given
by (\ref{zamena}). The Cartan involution of $\hgv$ can be
expressed as: \begin{equation} \omega (A\tp f({\bf{t}})+ \lambda_1
K_1 + \lambda_2 K_2 )=\tilde\omega(A) \tp f(w_0({\bf{t}}))-
\lambda_1 K_1 - \lambda_2 K_2.
\end{equation} It is not difficult to see that for any element
$g=A\tp f({\bf{t}}) \in L^{\slf}\g$ we have
 $(g, \, \omega(g)) >0$.

In QFT it is important to consider the correlations functions of
two or more fields. In 2-dimensional conformal field theory the
correlator of two fields \begin{equation}
<v^*|A(z)B(w)|v>\end{equation} is a rational function of $z$ and
$w$ in the domain $|z|>|w|$ with poles only at hyperplanes $z=0,\,
w=0,$ and $z=w.$ This important property in quantum field theory
in terms of OPE (operator product expansion) means that the
product of two fields at nearby points can be expanded in terms of
other fields and the small parameter $z-w.$ In terms of commutator
the OPE can be expressed as \begin{equation} A(z)B(w)=[A_- (z), \,
B(w)]+:\!A(z)B(w)\!:,\end{equation} where $[A_- (z), \, B(w)]$ is
a singular part of the OPE and \begin{equation} :\!A(z)B(w)\!:=
A(z)_+ B(w) +B(w)A(z)_-\end{equation} is the normal ordered
product of two fields. The singular part of the OPE of two local
fields $A(z)$ and $B(w),$ given by \begin{equation} [A_- (z),
\,B(w)]=\sum_{n=0}^{N} A_{(n)}B (w) \, \d^{(n)}_w \delta (z-w)_+
,\end{equation} is a finite linear combination of other fields
$A_{(n)}B (w)$ from the theory. Here $\d^{(n)}= \frac{1}{n!}
\,\d^n.$ The $\delta (z-w)_+$ is a regular function with respect
to $w$ in the domain $|z|>|w|$ and has the properties: $\d_w
\delta (z-w)_+=-\d_z \delta (z-w)_+$ and $\delta(z-w)_+
|_{w=0}=1/z.$ From these properties we have
\begin{equation}\label{sum} \delta (z-w)_+ = \frac{1}{z-w} \, \,
\,\,{\rm in \;\; the\;\; domain }\,\,\, |z|>|w|. \end{equation}
For the regular delta function we don't need these properties to
find the sum \begin{equation} \delta (z-w)_+ =\sum_{n=0}^{\infty}
z^{-n-1} w^n =\frac{1}{z-w}, \end{equation} but for the
generalized delta function this can give us some idea.

In order to have the consistent definition of the correlator of
two fields in the dimension higher then 2, we need to introduce
the OPE and the normal product of two fields on $M.$

Consider two fields $a^{\slf} ({\bf z}), \, b^{\slf} ({\bf w}) \in
\hgv.$ Define $a^{\slf} ({\bf z})_{\pm}$ as \begin{equation}
a^{\slf} ({\bf
 z})_{\pm}= a \tp \delta_{V}\, ({\bf z-w})_{\pm}, \end{equation} where
$\delta_{V}\, ({\bf z-w})_{\pm}$ are defined by (\ref{del+}). Then
the OPE can de defined as: \begin{equation}\label{opoe} a^{\slf
}({\bf z}) \, b^{\slf }({\bf w})= (a, \, b)\left
(\sum_{i=1}^{2}\,K_i \cdot
 L_{i}^{w}\, \delta_V\,({\bf z -w})_+
 \right )+
[a, \, b]^{\slf} ({\bf w})\, \delta_{V}\, ({\bf z -w})_{+} +\\
:\!a^{\slf} ({\bf z}) \, b^{\slf}({\bf w})\!: \end{equation} where
\begin{equation} :\!a^{\slf} ({\bf z}) \, b^{\slf} ({\bf w})\!:=
a^{\slf }({\bf z})_+ b^{\slf} ({\bf w}) +b^{\slf }({\bf
w})a^{\slf} ({\bf z})_- .\end{equation}

It would be nice to have the formula similar to (\ref{sum}) for
the $\delta_{V}\, ({\bf z-w})_{+},$ defined as (\ref{del+}). The
space $V^{\slf}$ is a sum of two dual subspaces $V^{\g} = V_+
\oplus V_-$. Here $\g=\slf $. The positive part of delta function
is defined as
\begin{equation} \delta_{V}\, ({\bf z- w})_{+}=\sum_{e_{\gamma}
\in V_+} e_{\gamma}^* \tp e_{\gamma} \in V_-^z \tp V_+^w
.\end{equation} The space $V_+$ of all regular functions on $M_+$
is the direct sum of two subspaces $V_+ = V_{++} \oplus V_{+-} \,$
where $V_{++}$ is a space of all regular functions on $U \simeq
\CC^3 .$ For the dual space $V_-$ we have the dual decomposition
$V_- = V_{-+} \oplus V_{--} .$ According to this decompositions we
can write \begin{equation} \delta_{V}\, ({\bf z- w})_{+} =
\delta_{V}\, ({\bf z- w})_{+ +} + \delta_{V}\, ({\bf z- w})_{+ -}.
\end{equation} The distribution \begin{equation} \delta_{V}\,({\bf
z- w})_{+ +} \in V_{--}^z \tp V_{++}^w \end{equation} is a regular
function with respect to ${\bf w}=(w_1, w_2, w_3)$ in the domain
$U \simeq \CC^3 .$ Thus, this function is well defined, when $w_i
=0. $ From the construction of the basis (\ref{ba}), we have
\begin{equation} \delta_{V}\,({\bf z- w})_{+ +} |_{w_1=0,
\,w_2=0,\, w_3 =0} =\phi (z_1, z_2, z_3) = \frac{1}{z_3 (z_3-z_1
z_2)}. \end{equation} Additionally, we have \begin{equation}
\d_i^w \, \delta_{V}\,({\bf z- w})_{+ +} = - \d_i^z \,
\delta_{V}\,({\bf z- w})_{+ +}, \,\,\,\, i=1, 2, 3.
\end{equation} This equality we understand at the level of distributions.

Consider the function of two sets of variables ${\bf z} =(z_1,
z_2, z_3)$ and ${\bf w}=(w_1, w_2, w_3)$ given by \begin{equation}
F({\bf z, \,w}) = \frac{z_1z_2-w_1w_2}{(z_1-w_1)(z_2-w_2)(z_3-w_3)
(z_3-z_1 z_2- (w_3-w_1 w_2))}=\end{equation}
$$\frac{z_1}{(z_1-w_1)(z_2-w_2) (z_3-z_1 z_2- (w_3-w_1 w_2))}
+\frac{w_2}{(z_2-w_2)(z_3-w_3) (z_3-z_1 z_2- (w_3-w_1 w_2))} .$$
\begin{proposition}
In the domain $D_1 :\;|z_1| >|w_1|>0,\;\;|z_2| >|w_2|>0,\;\;|z_3|
>|w_3|, \;\;\;\;|z_3-z_1 z_2|
> |w_3 -w_1 w_2 | $ the function $F({\bf z},{\bf w})=\delta_{V}\,({\bf z- w})_{+ +} $
in a sense that:\begin{equation} Res_{\bf z} \;F({\bf z},{\bf w})
\cdot f({\bf z}) = f({\bf z})_{++} \end{equation} and in the
domain $$D_2:\;|z_1| >|w_1|>0,\;\;|z_2| >|w_2|>0,\;\;|z_3| <|w_3|,
\,\, \,\,|z_3-z_1 z_2|
> |w_3 -w_1 w_2 | $$
the function $-F({\bf z},{\bf w})=\delta_{V}\,({\bf z- w})_{+ -} $
in a sense that:\begin{equation} Res_{\bf z} \;F({\bf z},{\bf w})
\cdot f({\bf z}) = -f({\bf z})_{+-} .\end{equation}
\end{proposition}

The proof is based on a realization of the trace as an integral
over the 3-dimensional cycle $S$ and the corresponding formal
power series expansion for $F({\bf z},{\bf w})$ on $S$
(\cite{G-K1}).

 Then, the OPE expansion of two fields $a({\bf z}), \, b({\bf
w}) \in \hgv$ can be written in terms of the function $F({\bf
z},{\bf w})$ if we replace in the OPE (\ref{opoe} ) the
$\delta_{V}\,({\bf z- w})_{+ }$ by the sum $$\;\;\;\;\;F({\bf
z},{\bf w})\;\;\; \;\;\;\;-\;\;\;\;\;\;\; F({\bf z},{\bf w})
$$
$$ \!\!\!\!\!\!\!\!\!\!\!\!\!\!\!\!\;\;\;\;|z_1| >|w_1|>0,\;\;|z_2| >|w_2|>0,
\;\;\;\;\;\;\;\;\;\;\;\;\;\
 |z_1| >|w_1|>0,\;\;|z_2| >|w_2|>0,$$
 $$\;\;\;|z_3|>|w_3|, \,\, \,\,|z_3-z_1 z_2|> |w_3 -w_1 w_2 |;\,\,\,\;\;\;|z_3|<|w_3|, \,\, \,\,|z_3-z_1 z_2|
> |w_3 -w_1 w_2 |$$
 This form of OPE is very similar to what we have in
2-dimensional conformal field theory.

\section{Generalized Virasoro algebra associated with the space
$V^{\slf }$}

In section 3 we constructed the $n$-products of the Virasoro
conformal algebra from the commutation relations of $\slt$
(\ref{vcon}). Now we will apply the same ideas to construct a
generalized Virasoro conformal algebra associated to the Lie
algebra $\slf .$

 For $\g =\slf$ positive
roots are $\{\al_1, \, \al_2,\,\, \al_3=\al_1 + \al_2 \}$ and the
corresponding weight basis in $\slf$ is $\{e_1, \, e_2,\, e_3, \,
h_1,\, h_2,\, f_1,\, f_2,\, f_3\}.$ Denote by $\d_i = ad \,(e_i).$
Lie algebra $\slf$ has a structure of the $U(\n_+)$-module with
the lowest weight $-\al_3.$ Via the ${\rm Poincar\acute
e}$-Birkhoff-Witt basis we can identify $U(\n_+) \simeq
\CC[\d_{1}, \d_2,\d_{3}].$ Denote by $L=f_3$ the lowest vector
with respect to this action, then $\slf \simeq \CC[\d_{1},
\d_2,\d_{3}] \,L$ is the rank one $\CC[\d_{1},
\d_2,\d_{3}]$-module. Consider a basis in $\slf$ of the form $$
L=f_3, \;\; \d_1 L= -f_2, \;\;\; \d_2 L= f_1; \;\;\;\; \d_1 \d_2
L= h_1.\;\;\; \d_3 L=h_1+h_2;$$\begin{equation} \d_1^{2} \d_2
L=-2e_1,\;\;\;\; \d_2 \d_3 L =-e_2,\;\;\;\;\d_{1} \d_2 \d_{3}L=
-e_3.\end{equation} The dual basis (with respect to the normalized
bilinear form ) is
\begin{equation}
 e_3,\;\;\; -e_2,\;\;\; e_1;\;\;\;\;\;
\frac{1}{3}h_1 - \frac{1}{3}h_2 , \;\;\;\; \frac{1}{3}h_1 +
\frac{2}{3}h_2;\;\;\;\;\; -\frac{1}{2}f_1,\;\;\;\; -f_2, \;\;\;\;
-f_3 .\end{equation}

For any element $e$ of the dual basis we define a bilinear
operation $L_{(e)}L$ product as: $$L_{(f_3^*)}L=[e_3,\,
L]=\d_3L,\;\;\; \;\;\;\;L_{((\d_1f_3)^*)}L=[-e_2,\,
L]=-\d_2L,$$\begin{equation}\label{confn}L_{((\d_2f_3)^*)}L=[e_1,\,
L]=\d_1 L,\;\;\;\;\;\;
 L_{((\d_3f_3)^*)}L=[-2h_1-h_2,
L]=3L.\end{equation} The other products are $0$ because $[f_i,\,
f_3]=0$ for $i=1,\,2,\,3.$

Let \begin{equation} \delta_{\scriptstyle{V }}
=\delta_{\scriptstyle{V^{\slf}}}({\bf t-z}) \end{equation}
 be a delta function associated with the space ${V^{\slf}}, $
 defined in (\ref{del1}). Denote by $L^\g ({\bf z})$ the field associated with the space
$V^{\slf} \simeq {\rm Fun }\, M$ with the commutation relation
defined by products (\ref{confn}) :$$ [L^\g({\bf z}),\,L^\g({\bf
w})]= \d_3 L^\g({\bf w})\,\delta_V({\bf z-w})-\d_2 L^\g({\bf
w})\,\d_1^w\delta_V({\bf z-w})+\d_1 L^\g({\bf
w})\,\d_2^w\delta_V({\bf z-w})+$$
\begin{equation}\label{vcr}
3 L^\g({\bf w})\,\d_3^w\delta_V({\bf z-w})\end{equation}
\begin{theorem}
  The defined above commutation relation satisfies to the skewsymmetry and
  Jacobi identity axioms
\end{theorem}

$ \bf Proof.$ We have $$[L^\g({\bf w}),\,L^\g({\bf z})]= \d_3
L^\g({\bf z})\,\delta_{{V}}({\bf z-w})-\d_2 L^\g({\bf
z})\,\d_1^z\delta_{V}({\bf z-w})+\d_1 L^\g({\bf
z})\,\d_2^z\delta_{V}({\bf z-w})+$$ $$3 L^\g({\bf
z})\,\d_3^z\delta_{V}({\bf z-w})= \d_3 L^\g({\bf
w})\,\delta_{V}({\bf z-w})+\d_2 L^\g({\bf
w})\,\d_1^w\delta_{V}({\bf z-w})+ $$ $$ \d_1\d_2 L^\g({\bf
w})\,\delta_{V}({\bf z-w})- \d_1 L^\g({\bf
w})\,\d_2^w\delta_{V}({\bf z-w})-\d_2\d_1 L^\g({\bf
w})\,\delta_{V}({\bf z-w})-$$ \begin{equation} 3 L^\g({\bf
w})\,\d_3^w\delta_{V}({\bf z-w})-3\,\d_3 L({\bf
w})\,\delta_{V}({\bf z-w})= -[L^\g({\bf z}),\,L^\g({\bf w})].
\end{equation}

Define a field on $V^{\slf}$ by \begin{equation}\label{fil} L^\g({\bf
z})=-\delta_{{V}}({\bf t-z})\,\d_3^t +\d_1^t\delta_{V}({\bf t-
z})\d_2^t-\d_2^t\delta_{V}({\bf t- z})\d_1^t .\end{equation} The
straightforward calculations show that this field satisfies to the
commutation relation (\ref{vcr}).

For any function $f\in V^{\slf}$ define an operator $ L_f$ as $
L^\g_f=Res_z \,(f({\bf z}) L^\g({\bf z}))$ so that
\begin{equation} L^\g_f =-f\,\d_3^t+\d_1f\,\d_2^t -\d_2 f\,\d_1^t
.\end{equation} From the direct calculations we have
\begin{equation}\label{brac} [L^\g_f,\,L^\g_g]=L^\g_{\{f,\,g\}},
\end{equation} where $\{f,\,g\}$ define a Lie bracket on the space $
V^{\slf}$ given by \begin{equation}\label{brack} \{f,\,g\}=
g\,\d_3f-f\,\d_3g +\d_1f\, \d_2g -\d_1g\,\d_2f .\end{equation}
This bracket is skewsymmetric and satisfies the Jacobi identity
\begin{equation} \{h , \,\{f,\,g\}\} + \{f , \,\{g,\,h\}\}+\{g ,
\,\{h,\,f\}\}=0\end{equation} for any $h,\,f,\,g \in V^\g,$ so it
is a Lie bracket. This completes the proof.

\begin{proposition} The Lie bracket defined in (\ref{brack}) corresponds
to the contact Lie bracket on the real space
$M_{\RR}=V^{\g_{\RR}}$( we now suppose that $(a,\,b,\, c)$ are
real coordinates).  The operators $L^\g_f$ are contact vector
fields on $M_{\RR}$ with the contact Hamiltonian $f.$ The
corresponding contact structure is defined by a 1-form $\beta ,$
such that $ \beta (\d_3)=-1,\;\;\;\beta (\d_1)=0,\;\;\;\beta
(\d_2)=0,\;\;\; d\beta(\d_3,\, \cdot)=0$ and $\beta \wedge
d\,\beta \ne 0.$
\end{proposition}
{\bf Proof.} Let $\beta$ be a contact 1-form on 3-dimensional
 real manifold $M :$  $\,\beta \wedge d\,\beta \ne 0.$
Vector field $\xi_f$ on $M$ is contact with the contact
hamiltonian $f\in {\rm Fun}\, M,$ if $L_{\xi_f}$ preserves the
foliation $\beta =0$ and $\beta(\xi_f)=f$. Then $\xi_f=f\cdot
\al+\theta(f)$ where $\al$ is the vector field defined by
$\beta(\al)=1$ and $d\beta(\al,\,\cdot)=0$ and $\theta$ is the
bivector field on $M,$ inverse  to the 2-form $d\beta .$ The
contact Lie bracket on ${\rm Fun}\,M$ is defined by
$[\xi_f,\,\xi_g]= \xi_{\{f,\,g\}}$ or more explicitly:
\begin{equation}
\{f,\,g\}=f\al(g)-g\al(f) +\theta\,(f,\,g).
\end{equation}
Take $\al=-\d_3, \;\;\;\theta=\d_1\wedge \d_2$ and
$\beta=-dc+adb$. This completes the proof.
\begin{proposition} The Lie algebra $V^{\slf}$ with the bracket (\ref{brack})
contains the Lie subalgebra isomorphic to $\slf .$
\end{proposition}
$ \bf Proof.$ For simplicity we will use the homogeneous
coordinates $(a,\,b,\,c)$ on  $M^{\g}$ The action of $\slf$ in
terms of these coordinates is given by (\ref{action}). Operators
\begin{equation} L_1,\;\;\;L_a,\;\;\; L_b ,\;\;\; L_c,
\;\;\;L_{c-ab},\;\;\;L_{ac},\;\;\;L_{-b(c-ab)},\;\;\;
L_{c(c-ab)}\end{equation} are closed with respect to the Lie
bracket (\ref{brac}) and form the Lie algebra isomorphic to the
$\slf.$ The Cartan subalgebra generators $L_c = -c\, \d_3 - a \,
\d_1, \,\,\,\,\,\, L_{c-ab}= -(c-ab)\, \d_3 -b\, \d_2 $ are
grading operators from the previous section and they are analogues
of the energy operator $L_0$ in two dimensional CFT.

 We constructed a natural higher dimensional
analogues of the Virasoro algebra and the Virasoro conformal
algebra associated with the space $V^{\slf}.$ We suggest to denote
this algebra by $Vir^{\slf} .$

 We also define a semidirect sum of $Vir^{\slf} $ and the generalized
 affine Kac-Moody algebra $\hgv$ associated with
the space $V^{\slf}.$ The $Vir^{\slf} $ acts on $\hgv $ via its
action as derivations of $V^{\slf} .$ In terms of fields this
action is given by $$ [L^\g({\bf z}),\, a^{\slf}({\bf w})]=\d_3
a^{\slf}({\bf w})\,\delta_{{V}}({\bf z-w})-\d_2 a^{\slf}({\bf
w})\,\d_1^w\delta_{V}({\bf z-w})+\d_1 a^{\slf}({\bf
w})\,\d_2^w\delta_{V}({\bf z-w})+$$\begin{equation} 2a^{\slf}({\bf
w})\,\d_3^w\delta_{V}({\bf z-w}) \end{equation} where $
a^{\slf}({\bf w}) \in \hgv .$

\end{document}